\documentstyle[amssymb,12pt]{article}

\begin{document}

\title{Quantum deformation of Whittaker modules and Toda lattice}

\author{Alexey Sevostyanov \footnote{e-mail seva@teorfys.uu.se}\\ 
Institute of Theoretical Physics, Uppsala University
}

\maketitle

\begin{abstract}
In 1978 Kostant suggested the Whittaker model
of the center of the universal enveloping algebra $U({\frak g})$
of a complex simple  Lie algebra $\frak g$. 
The main result is that the center of $U({\frak g})$ is isomorphic
to a commutative subalgebra in $U({\frak b}_+)$, where 
${\frak b}_+$ is a Borel subalgebra in $\frak g$. 
This observation is used in the theory of principal series
representations of the corresponding Lie group $G$ and in 
the proof of complete integrability of the quantum Toda lattice.
 In this paper we generalize the Kostant's construction
to quantum groups. In our construction we use quantum analogues of
regular nilpotent elements defined in \cite{S1}. 
 Using the Whittaker model of the center of 5the algebra $U_h({\frak g})$ 
we define quantum deformations of Whittaker modules.
The new Whittaker model is also applied to the deformed quantum Toda
lattice recently studied by Etingof in \cite{Et}. We give new proofs of his 
results
which resemble the original Kostant's proofs for the quantum Toda lattice.
\end{abstract}

\renewcommand{\theequation}{\thesection.\arabic{equation}}

\newtheorem{theorem}{Theorem}{}
\newtheorem{lemma}[theorem]{Lemma}{}
\newtheorem{corollary}[theorem]{Corollary}{}
\newtheorem{conjecture}[theorem]{Conjecture}{}
\newtheorem{proposition}[theorem]{Proposition}{}
\newtheorem{axiom}{Axiom}{}
\newtheorem{remark}{Remark}{}
\newtheorem{example}{Example}{}
\newtheorem{exercise}{Exercise}{}
\newtheorem{definition}{Definition}{}

\renewcommand{\thetheorem}{\thesection.\arabic{theorem}}

\renewcommand{\thelemma}{\thesection.\arabic{lemma}}

\renewcommand{\theproposition}{\thesection.\arabic{proposition}}

\renewcommand{\thecorollary}{\thesection.\arabic{corollary}}

\renewcommand{\theremark}{\thesection.\arabic{remark}}

\renewcommand{\thedefinition}{\thesection.\arabic{definition}}

\setcounter{equation}{0}
\setcounter{theorem}{0}

\section*{Introduction}

In 1978 Kostant suggested the {\em Whittaker model}
of the center of the universal enveloping algebra $U({\frak g})$
of a complex simple  Lie algebra ${\frak g}$. An essential role
in this construction is played by a non--singular character $\chi$
of the maximal nilpotent subalgebra ${{\frak n}_+} \subset \frak{g}$.
The main result is that the center of $U({\frak g})$ is isomorphic
to a commutative subalgebra in $U({\frak b}_-)$, where 
${\frak b}_- \subset {\frak g}$ is the opposite Borel subalgebra. 
This observation is used in the theory of principal series
representations of the corresponding Lie group $G$ and in 
the proof of complete integrability of the quantum Toda lattice.

The goal of this paper is to generalize the Kostant's construction
to quantum groups. An obvious obstruction is the fact that the
subalgebra in $U_h({\frak g})$ 
generated by positive root generators (subject to the quantum Serre relations)
does not have non-singular characters. In order to overcome
this difficulty we use a family of new realizations of
quantum groups introduced in \cite{S1}. 
The modified quantum Serre relations allow for non--singular
characters, and we are able to construct the Whittaker model
of the center of $U_h({\frak g})$.

Using the Whittaker model of the center of $U_h({\frak g})$ we introduce
quantum deformations of Whittaker modules.

The new Whittaker model is also applied to the deformed quantum Toda
lattice recently studied by Etingof (see \cite{Et}). We give new proofs of his 
results
which resemble the original Kostant's proofs for the quantum Toda lattice.

The paper is organized as follows. Section \ref{Witt} contains a review of
Kostant's results on the Whittaker model and Whittaker modules \cite{K}, 
\cite{K1}. 
In order to create a pattern for 
proofs in the  quantum group case we recall most of the Kostant's proofs.
The central part of the paper is Section \ref{qWitt}. There we discuss 
properties of
new realizations of finite-dimensional quantum groups and present
the Whittaker model of the center of $U_h({\frak g})$. 
In Section
\ref{Whittmodh} we introduce quantum deformed Whittaker modules. Section 
\ref{toda}
contains a discussion of the deformed quantum Toda lattice.


\section{Whittaker model}\label{Witt}

In this section we recall the Whittaker model of the center
of the universal enveloping algebra $U({\frak g})$, where ${\frak g}$ is a 
complex simple Lie algebra.


\subsection{Notation}\label{notation}

Fix the notation used throughout of the text.
Let $G$ be a
connected simply connected finite--dimensional complex simple Lie group, $%
{\frak g}$ its Lie algebra. Fix a Cartan subalgebra ${\frak h}\subset {\frak %
g}\ $and let $\Delta $ be the set of roots of $\left( {\frak g},{\frak h}%
\right) .$ Choose an ordering in the root system. Let $\alpha_i,~i=1,\ldots 
l,~~l=rank({\frak g})$ be the 
simple roots, $\Delta_+=\{ \beta_1, \ldots ,\beta_N \}$  
the set of positive roots. Denote by $\rho$ a half of the sum of positive roots,
$\rho=\frac 12 \sum_{i=1}^N\beta_i$.
Let $H_1,\ldots ,H_l$ be the set of simple root generators of $\frak h$. 

Let $a_{ij}$ be the corresponding Cartan matrix.
Let $d_1,\ldots , d_l$ be coprime positive integers such that the matrix 
$b_{ij}=d_ia_{ij}$ is symmetric. There exists a unique non--degenerate invariant
symmetric bilinear form $\left( ,\right) $ on ${\frak g}$ such that 
$(H_i , H_j)=d_j^{-1}a_{ij}$. It induces an isomorphism of vector spaces 
${\frak h}\simeq {\frak h}^*$ under which $\alpha_i \in {\frak h}^*$ corresponds 
to $d_iH_i \in {\frak h}$. We denote by $\alpha^\vee$ the element of $\frak h$ 
that 
corresponds to $\alpha \in {\frak h}^*$ under this isomorphism.
The induced bilinear form on ${\frak h}^*$ is given by
$(\alpha_i , \alpha_j)=b_{ij}$.

Let $W$ be the Weyl group of the root system $\Delta$. $W$ is the subgroup of 
$GL({\frak h})$ 
generated by the fundamental reflections $s_1,\ldots ,s_l$,
$$
s_i(h)=h-\alpha_i(h)H_i,~~h\in{\frak h}.
$$
The action of $W$ preserves the bilinear form $(,)$ on $\frak h$. 
We denote a representative of $w\in W$ in $G$ by
the same letter. For $w\in W, g\in G$ we write $w(g)=wgw^{-1}$.

Let ${{\frak b}_+}$ be the positive Borel subalgebra and ${\frak b}_-$
the opposite Borel subalgebra; let ${\frak n}_+=[{{\frak b}_+},{{\frak b}_+}]$ 
and $%
{\frak n}_-=[{\frak b}_-,{\frak b}_-]$ be their 
nil-radicals. Let $H=\exp {\frak h},N_+=\exp {{\frak n}_+},
N_-=\exp {\frak n}_-,B_+=HN_+,B_-=HN_-$ be
the Cartan subgroup, the maximal unipotent subgroups and the Borel subgroups
of $G$ which correspond to the Lie subalgebras ${\frak h},{{\frak n}_+},%
{\frak n}_-,{\frak b}_+$ and ${\frak b}_-,$ respectively.

We identify $\frak g$ and its dual by means of the canonical invariant bilinear 
form. 
Then the coadjoint 
action of $G$ on ${\frak g}^*$ is naturally identified with the adjoint one. We 
also identify 
${{\frak n}_+}^*\cong {\frak n}_-,~{{\frak b}_+}^*\cong {\frak b}_-$. 

Let ${\frak g}_\beta$ be the root subspace corresponding to a root $\beta \in 
\Delta$, 
${\frak g}_\beta=\{ x\in {\frak g}| [h,x]=\beta(h)x \mbox{ for every }h\in 
{\frak h}\}$.
${\frak g}_\beta\subset {\frak g}$ is a one--dimensional subspace. 
It is well--known that for $\alpha\neq -\beta$ the root subspaces ${\frak 
g}_\alpha$ and ${\frak g}_\beta$ are orthogonal with respect 
to the canonical invariant bilinear form. Moreover ${\frak g}_\alpha$ and 
${\frak g}_{-\alpha}$
are non--degenerately paired by this form.

Root vectors $X_{\alpha}\in {\frak g}_\alpha$ satisfy the following relations:
$$
[X_\alpha,X_{-\alpha}]=(X_\alpha,X_{-\alpha})\alpha^\vee.
$$


\subsection{The Whittaker model}\label{whitt}

\setcounter{equation}{0}
\setcounter{theorem}{0}

In this section we introduce the Whittaker model of the center of the universal 
enveloping
algebra $U({\frak g})$. We start by recalling the classical result of Chevalley
which describes the structure of the center.

Let $Z({\frak g})$ be the center of $U({\frak g})$. The standard filtration 
$U_k({\frak g})$ 
in $U({\frak g})$ induces a filtration $Z_k({\frak g})$ in $Z({\frak g})$. The 
following
important theorem may be found for instance in \cite{Bur1}, Ch.8, \S 8, no. 3, 
Corollary 1 and no.5, Theorem 2.
\vskip 0.3cm
\noindent
{\bf Theorem (Chevalley)}
{\em One can choose
elements $I_k\in Z_{m_k+1}({\frak g}),~~k=1,\ldots l$, where $m_k$ are called 
the exponents of $\frak g$,
such that
$Z({\frak g})={\Bbb C}[I_1,\ldots , I_l]$ is a polynomial algebra in $l$ 
generators.}
\vskip 0.3cm
The adjoint action of $G$ on $\frak g$ naturally extends to $S({\frak g})$. 
Let $S({\frak g})^G$ be the algebra of $G$--invariants in $S({\frak g})$.
Clearly, $GrZ({\frak g})\cong S({\frak g})^G$. In particular $S({\frak 
g})^G\cong 
{\Bbb C}[\widehat I_1,\ldots , \widehat I_l]$, where $\widehat I_i=Gr 
I_i,~i=1,\ldots ,l$. The elements
$\widehat I_i,~i=1,\ldots ,l$ are called fundamental invariants.

Following Kostant we shall realize the center $Z({\frak g})$ of the universal 
enveloping 
algebra $U({\frak g})$ as a subalgebra in $U({\frak b}_-)$.
Let 
$$
\chi :{{\frak n}_+} \rightarrow {\Bbb C}
$$
be a character of ${{\frak n}_+}$. 
Since ${{\frak n}_+}=\sum_{i=1}^l{\Bbb C}X_{\alpha_i} \oplus[{{\frak 
n}_+},{{\frak n}_+}]$
it is clear that $\chi$ is completely determined by the constants $c_i=\chi 
(X_{\alpha_i}),~i=1,\ldots ,l$ and 
$c_i$ are arbitrary. In \cite{K} $\chi$ is called non--singular if $c_i\neq 0$ 
for all $i$.

Let $f=\sum_{i=1}^l X_{-\alpha_i}\in {\frak n}_-$ be a regular nilpotent 
element.
From the properties of the invariant bilinear form (see Section \ref{notation}) 
it follows that
$(f,[{{\frak n}_+},{{\frak 
n}_+}])=0,~~(f,X_{\alpha_i})=(X_{-\alpha_i},X_{\alpha_i})$, 
and hence the map
$x\mapsto (f,x),~~x\in {{\frak n}_+}$ is a non--singular character of ${{\frak 
n}_+}$.

Recall that in our choice of root vectors no normalization was made. But now 
given a non--singular
character $\chi :{{\frak n}_+}\rightarrow {\Bbb C}$ we will say that $f$ 
corresponds to $\chi$ in case
$$
\chi (X_{\alpha_i}) =(X_{-\alpha_i},X_{\alpha_i}).
$$
Conversely if $\chi$ is non--singular there is a unique choice of $f$ so that 
$f$ corresponds to $\chi$. 
In this case $\chi (x)=(f,x)$ for every $x\in {{\frak n}_+}$.

Naturally, the character $\chi$ extends to a character of 
the universal enveloping algebra $U({{\frak n}_+})$. 
Let $U_\chi ({{\frak n}_+})$ be the kernel of this extension so that 
one has a direct sum
$$
U({{\frak n}_+})={\Bbb C}\oplus U_\chi ({{\frak n}_+}).
$$

 Since ${\frak g}={\frak b}_-\oplus {{\frak n}_+}$ we have a linear 
isomorphism $U({\frak g})=U({\frak b}_-)\otimes U({{\frak n}_+})$ and hence 
the direct sum
\begin{equation}\label{maindec}
U({\frak g})=U({\frak b}_-) \oplus I_\chi,
\end{equation}
where $I_\chi=U({\frak g})U_\chi ({{\frak n}_+})$ is the left--sided ideal 
generated by 
$U_\chi ({{\frak n}_+})$.

For any $u\in U({\frak g})$ let $u^\chi\in U({\frak b}_-)$ be its component in 
$U({\frak b}_-)$ relative to the decomposition (\ref{maindec}). Denote by 
$\rho_\chi$
the linear map 
$$
\rho_\chi : U({\frak g}) \rightarrow U({\frak b}_-)
$$
given by $\rho_\chi (u)=u^\chi$.
Let $W({\frak b}_-)=\rho_\chi (Z({\frak g}))$. 
\vskip 0.3cm
\noindent
{\bf Theorem A (\cite{K}, Theorem 2.4.2)}
{\em The map
\begin{equation}\label{map}
\rho_\chi : Z({\frak g}) \rightarrow W({\frak b}_-)
\end{equation}
is an isomorphism of algebras. In particular 
$$
W({\frak b}_-)={\Bbb C}[I_1^\chi ,\ldots , I_l^\chi ], 
~~I_i^\chi=\rho_\chi(I_i),~~i=1,\ldots ,l
$$
is a polynomial algebra in $l$ generators.}
\vskip 0.3cm
\noindent
{\em Proof.}
First, we show that the map (\ref{map}) is an algebra homomorphism.
If $u,v\in Z({\frak g})$ then $u^\chi v^\chi \in U({\frak b}_-)$ and
$$
uv-u^\chi v^\chi =(u-u^\chi )v+u^\chi (v-v^\chi ).
$$
Since $(u-u^\chi )v=v(u-u^\chi )$ the r.h.s. of the last equality is an element 
of $I_\chi$.
This proves $u^\chi v^\chi =(uv)^\chi$.

By definition the map (\ref{map}) is surjective. We have to prove that it is 
injective.
Let $U({\frak g})^{\frak h}$ be the centralizer of $\frak h$ in $U({\frak g})$. 
Clearly 
$Z({\frak g})\subseteq U({\frak g})^{\frak h}$. From the 
Poincar\'{e}--Birkhoff--Witt 
theorem it follows that every element $z\in U({\frak g})^{\frak h}$ may be 
uniquely
written as
$$
z=\sum_{p,q\in{\Bbb N}^N,<p>=<q>}X_{-\beta_1}^{p_1}\ldots 
X_{-\beta_N}^{p_N}\varphi_{p,q}
X_{\beta_1}^{q_1}\ldots X_{\beta_N}^{q_N},
$$
where $<p>=\sum_{i=1}^r p_i \beta_i \in {\frak h}^*$ and $\varphi_{p,q} \in 
U({\frak h})$. 

Now recall that $\chi (X_{\beta_i})=0$ if $\beta_i$ is not a simple root, and we 
easily obtain
$$
\rho_\chi (z)=\sum_{p,q\in{\Bbb N}^l,<p>=<q>\neq 
0}X_{-\alpha_{k_1}}^{p_{j_1}}\ldots X_{-\alpha_{k_l}}^{p_{j_l}}\varphi_{p,q}
\prod_{i=1}^lc_{k_i}^{q_{j_i}}+\varphi_{0,0}.
$$

Let $z\in Z({\frak g})$. One knows that the map 
$$
Z({\frak g})\rightarrow U({\frak h}),~~z\mapsto \varphi_{0,0},
$$
called the Harich-Chandra homomorphism, is injective (see (c), p. 232 in 
\cite{Dix}). It follows that 
the map (\ref{map}) is also injective.
\begin{remark}\label{homrho}
The first part of the proof of Theorem A only used the fact that $v\in Z({\frak 
g})$. Therefore 
$$
\rho_\chi(uv)=\rho_\chi(u)\rho_\chi(v)
$$
for any $u\in U({\frak g}),~~v\in Z({\frak g})$.
\end{remark} 
\vskip 0.3cm
\noindent
{\bf Definition A}
{\em The algebra $W({\frak b}_-)$ is called the Whittaker model of $Z({\frak 
g})$.}  
\vskip 0.3cm

Next we equip $U({\frak b}_-)$ with a structure of a left $U({{\frak n}_+})$ 
module in such a
way that $W({\frak b}_-)$ is realized as the space of invariants with respect to 
this action.

Let $Y_\chi$ be the left $U({\frak g})$ module defined by 
$$
Y_\chi =U({\frak g})\otimes_{U({{\frak n}_+})}{\Bbb C}_\chi ,
$$
where ${\Bbb C}_\chi$ denotes the 1--dimensional $U({{\frak n}_+})$--module 
defined by $\chi$.
Obviously $Y_\chi$ is just the quotient module $U({\frak g})/I_\chi$. From 
(\ref{maindec}) it 
follows that the map
\begin{equation}\label{iso2}
U({\frak b}_-)\rightarrow Y_\chi;~~v\mapsto v\otimes 1
\end{equation}
is a linear isomorphism. 

It is convenient
to carry the module structure of $Y_\chi$ to $U({\frak b}_-)$. For $u\in 
U({\frak g}),~~
v\in U({\frak b}_-)$ the induced action $u\circ v$ has the form
\begin{equation}\label{indact}
u\circ v=(uv)^\chi.
\end{equation}
The restriction of this action to $U({{\frak n}_+})$ may be changed by tensoring 
with 1--dimensional
$U({{\frak n}_+})$--module defined by $-\chi$. That is $U({\frak b}_-)$ becomes 
an $U({{\frak n}_+})$
module where if $x\in {U({\frak n}_+)},~~v\in U({\frak b}_-)$ one puts
\begin{equation}\label{mainact}
x\cdot v=x\circ v-\chi (x)v.
\end{equation}
\vskip 0.3cm
\noindent
{\bf Lemma A (\cite{K}, Lemma 2.6.1.)}
{\em Let $v\in U({\frak b}_-)$ and $x\in {U({\frak n}_+)}$. Then}
$$
x\cdot v =[x,v]^\chi.
$$
\vskip 0.3cm
\noindent
{\em Proof.} 
By definition $x\cdot v=(xv)^\chi -\chi (x)v$. Then we have $xv=[x,v]+vx$ and 
hence
$x\cdot v=([x,v])^\chi +(vx)^\chi -\chi (x)v$. But clearly $(vx)^\chi =v\chi 
(x)$. Thus
$x\cdot v=([x,v])^\chi$.

The action (\ref{mainact}) may be lifted to an action of the unipotent group 
$N_+$.
Consider the space $U({\frak b}_-)^{N_+}$ of $N_+$ invariants 
in $U({\frak b}_-)$ with respect to this
action. Clearly, 
$W({\frak b}_-)\subseteq U({\frak b}_-)^{N_+}$ . 
\vskip 0.3cm
\noindent
{\bf Theorem B (\cite{K}, Theorems 2.4.1, 2.6)}
{\em Suppose that the character $\chi$ is non--singular. Then the space of $N_+$ 
invariants 
in $U({\frak b}_-)$ with respect to the
action (\ref{mainact}) is isomorphic to $W({\frak b}_-)$, i.e.}
\begin{equation}\label{inv}
U({\frak b}_-)^{N_+}\cong W({\frak b}_-).
\end{equation}
\vskip 0.3cm

%

\subsection{Whittaker modules}\label{Whittmod}

\setcounter{equation}{0}
\setcounter{theorem}{0}

In this section we recall basic facts on Whittaker modules (see \cite{K}).

Let $V$ be a $U({\frak g})$ module. The action is denoted by $uv$ for $u\in 
U({\frak g}),~~v\in V$. Let 
$\chi:{\frak n}_+ \rightarrow {\Bbb C}$ be a non--singular character of ${\frak 
n}_+$ (see Section \ref{whitt}).
A vector $w\in V$ is called a Whittaker vector (with respect to $\chi$) if
$$
xw=\chi (x)w
$$
for all $x\in U({\frak n}_+)$. A Whittaker vector $w$ is called a cyclic 
Whittaker vector (for $V$) if 
$U({\frak g})w=V$. A $U({\frak g})$ module $V$ is called a Whittaker module if 
it contains a cyclic 
Whittaker vector.

If $V$ is any $U({\frak g})$ module we let $U_V({\frak g})$ be the annihilator 
of $V$. Then $U_V({\frak g})$
defines a central ideal $Z_V({\frak g})$ by putting
\begin{equation}\label{zv}
Z_V({\frak g})=Z({\frak g})\cap U_V({\frak g}).
\end{equation}

Now assume that $V$ is a Whittaker module for $U({\frak g})$ and $w\in V$ is a 
cyclic Whittaker vector.
Let $U_w({\frak g})\subseteq U({\frak g})$ be the annihilator of $w$. Thus 
$U_V({\frak g}) \subseteq
U_w({\frak g})$, where $U_w({\frak g})$ is a left ideal and $U_V({\frak g})$ is 
a two--sided ideal in 
$U({\frak g})$. One has $V=U({\frak g})/U_w({\frak g})$ as $U({\frak g})$ 
modules so that $V$ is determined
up to equivalence by $U_w({\frak g})$. Clearly $I_\chi=U({\frak g})U_\chi 
({{\frak n}_+})\subseteq 
U_w({\frak g})$ and $U({\frak g})Z_V({\frak g})\subseteq U_w({\frak g})$.

The following theorem says that up to equivalence $V$ is determined by the 
central ideal $Z_V({\frak g})$.
\vskip 0.3cm
\noindent
{\bf Theorem F (\cite{K}, Theorem 3.1)}
{\em Let $V$ be any $U({\frak g})$ module which admits a cyclic Whittaker vector 
$w$ and let
$U_w({\frak g})$ be the annihilator of $w$. Then}
\begin{equation}\label{anndec}
U_w({\frak g})=U({\frak g})Z_V({\frak g})+I_\chi .
\end{equation}
\vskip 0.3cm

The proof of Theorem F is based on the following lemma. We use the notation of 
Section \ref{whitt}.
If $X\subseteq U({\frak g})$ let $X^\chi =\rho_\chi(X)$. Note that $U_w({\frak 
g})$ is stable under the map
$u\mapsto \rho_\chi(u)$. We recall also that by Theorem $\rm A$ $\rho_\chi$ 
induces an algebra isomorphism
$Z({\frak g}) \rightarrow W({\frak b}_-)$, where $W({\frak b}_-)=Z({\frak 
g})^\chi$. Thus if $Z_*$ is any ideal in
$Z({\frak g})$ then $W_*({\frak b}_-)=Z_*^\chi$ is an isomorphic ideal in 
$W({\frak b}_-)$.
But $(U({\frak g})Z_*)^\chi=U({\frak b}_-)W_*({\frak b}_-)$ by Remark 
\ref{homrho}. Thus by (\ref{maindec})
one has the direct sum
\begin{equation}\label{anndec*}
U({\frak g})Z_*+I_\chi=U({\frak b}_-)W_*({\frak b}_-)\oplus I_\chi.
\end{equation}
\vskip 0.3cm
\noindent
{\bf Lemma B (\cite{K}, Lemma 3.1)}
{\em Let $X=\{ v \in U({\frak b}_-) | (x\cdot v)w=0 \mbox{ for all }x \in {\frak 
n}_+\}$, where
$x\cdot v$ is given by (\ref{mainact}). Then
\begin{equation}
X=U({\frak b}_-)W_V({\frak b}_-)+W({\frak b}_-),
\end{equation}
where $W_V({\frak b}_-)=Z_V({\frak g})^\chi$. Furthermore if we denote 
$U_w({\frak b}_-)=U_w({\frak g})\cap U({\frak b}_-)$ then
\begin{equation}\label{bann}
U_w({\frak b}_-)=U({\frak b}_-)W_V({\frak b}_-).
\end{equation}
}
\vskip 0.3cm
\noindent
{\em Proof of Theorem F.} Let $u\in U_w({\frak g})$. We wish to show that 
$u\in U({\frak g})Z_V({\frak g})+I_\chi$. But by (\ref{anndec*}) it suffices to 
show that
$u^\chi \in U({\frak b}_-)W_V({\frak b}_-)$. Since $u^\chi \in U_w({\frak b}_-)$ 
the result 
follows from (\ref{bann}).
\vskip 0.3cm

Now one can determine, up to equivalence, the set of all Whittaker modules. They 
are naturally
parametrized by the set of all ideals in the center $Z({\frak g})$.
\vskip 0.3cm
\noindent
{\bf Theorem G (\cite{K}, Theorem 3.2)} 
{\em Let $V$ be any Whittaker module for $U({\frak g})$, the universal 
enveloping algebra of a 
simple Lie algebra $\frak g$. Let $U_V({\frak g})$ be the annihilator of $V$ and 
let $Z({\frak g})$
be the center of $U({\frak g})$. Then the correspondence
\begin{equation}\label{idcorr}
V\mapsto Z_V({\frak g}),
\end{equation}
where $Z_V({\frak g})=U_V({\frak g})\cap Z({\frak g})$, sets up a bijection 
between the set of all 
equivalence classes of Whittaker modules and the set of all ideals in $Z({\frak 
g})$.}
\vskip 0.3cm
\noindent
{\em Proof.} Let $V_i,~~i=1,2$ be two Whittaker modules. If $Z_{V_1}({\frak 
g})=Z_{V_2}({\frak g})$ then
clearly $V_1$ is equivalent to $V_2$ by (\ref{anndec}). Thus the map 
(\ref{idcorr}) is injective on
equivalence classes. 

Conversely, let $Z_*$ be any ideal in $Z({\frak g})$ and let $L=U({\frak 
g})Z_*+I_\chi$. Then
$V=U({\frak g})/L$ is a Whittaker module, where we can take $U_w({\frak g})=L$. 
But then 
$L=U({\frak g})Z_V({\frak g})+I_\chi$ by Theorem F. By (\ref{anndec}) this 
implies $Z_V({\frak g})^\chi=
Z_*^\chi$. However by Theorem A $\rho_\chi$ is injective on $Z({\frak g})$. 
Therefore $Z_V({\frak g})=
Z_*$. Hence the map (\ref{idcorr}) is surjective.

Now consider the subalgebra $Z({\frak g})U({\frak n}_+)$ in $U({\frak g})$. By 
Theorem 2.1 in \cite{K1}
$Z({\frak g})U({\frak n}_+)\cong Z({\frak g})\otimes U({\frak n}_+)$. Now let 
$Z_*$ be any ideal in
$Z({\frak g})$ and regard $Z({\frak g})/Z_*$ as a $Z({\frak g})$ module. Equip 
$Z({\frak g})/Z_*$ with a 
structure of $Z({\frak g})\otimes U({\frak n}_+)$ module by $u\otimes v 
y=\chi(v)uy$, where $u\in Z({\frak g}),~~
v\in U({\frak n}_+),~~y\in Z({\frak g})/Z_*$. We denote this module by 
$(Z({\frak g})/Z_*)_\chi$

The following result is another way of expressing Theorem G. 
\vskip 0.3cm
\noindent
{\bf Theorem H (\cite{K}, Theorem 3.3)}
{\em Let $V$ be any $U({\frak g})$ module. Then $V$ is a Whittaker module if and 
only if one has an isomorphism
\begin{equation}\label{whittind}
V\cong U({\frak g})\otimes_{Z({\frak g})\otimes U({\frak n}_+)}(Z({\frak 
g})/Z_*)_\chi
\end{equation}
of $U({\frak g})$ modules. Furthermore in such a case the ideal $Z_*$ is unique 
and is given by 
$Z_*=Z_V({\frak g})$, where $Z_V({\frak g})$ is defined by (\ref{zv}).}
\vskip 0.3cm
\noindent
{\em Proof.}If $1_*$ is the image of $1$ in $Z({\frak g})/Z_*$ then the 
annihilator in 
$Z({\frak g})\otimes U({\frak n}_+)$ of $1_*$ is $U({\frak n}_+)Z_*+Z({\frak 
g})U_\chi ({{\frak n}_+})$.
Thus the annihilator in $U({\frak g})$ of $1\otimes 1_*=w$ in the right side of 
(\ref{whittind}) is
$U({\frak g})Z_*+I_\chi$. The result then follows from Theorem G since $w$ is 
clearly a cyclic
generator of this module.

Now one can determine all the Whittaker vectors in a Whittaker module and 
duscuss the question of
irreducibility for Whittaker modules.
\vskip 0.3cm
\noindent
{\bf Theorem K (\cite{K}, Theorem 3.4)}
{\em Let $V$ be any $U({\frak g})$ module with a cyclic Whittaker vector $w\in 
V$. Then any $v\in V$ is
a Whittaker vector if and only if $v$ is of the form $v=uw$, where $u\in 
Z({\frak g})$. Thus the 
space of all Whittaker vectors in $V$ is a cyclic $Z({\frak g})$ module which is 
isomorphic to 
$Z({\frak g})/Z_V({\frak g})$.}
\vskip 0.3cm
\noindent
{\em Proof.} Obviously if $v=uw$ for $u\in Z({\frak g})$ then $v$ is a Whittaker 
vector. Conversely let
$v\in V$ be a Whittaker vector. Write $v=uw$ for $u\in U({\frak g})$. Then 
clearly $v=u^\chi w$ so that
we can assume $u\in U({\frak b}_-)$. But now if $x\in {\frak n}_+$ then 
$xuw=\chi(x)uw$. But also
$usw=\chi(x)uw$. Thus $[x,u]w=0$ and hence $[x,u]^\chi w=0$. But $x\cdot 
u=[x,u]^\chi$ by Lemma A. Thus
in the notation of Lemma B one has $u\in X$. But then by Lemma B one can write 
$u=u_1+u_2$, where 
$u_1\in U_w({\frak b}_-)$ and $u_2 \in W({\frak b}_-)$. But then $u_1w=0$. Thus 
$v=u_2w$. But now $u_2=u_3^\chi$,
where $u_3\in Z({\frak g})$ by Theorem A. But then $v=u_3w$ which proves the 
theorem.

If $V$ is any $U({\frak g})$ module then ${\rm End}_{U}~V$ denotes the algebra 
of operators on $V$ which
commute with the action of $U({\frak g})$. If $\pi_V: U({\frak g})\rightarrow 
{\rm End}~V$ is the representation 
defining the $U({\frak g})$ module structure on $V$ then clearly $\pi_V(Z({\frak 
g}))\subseteq {\rm End}_{U}~V$.
Furthermore it is also clear that $\pi_V(Z({\frak g}))\cong Z({\frak 
g})/Z_V({\frak g})$.
\vskip 0.3cm
\noindent
{\bf Theorem L (\cite{K}, Theorem 3.5)}
{\em Assume that $V$ is a Whittaker module. Then ${\rm End}_{U}~V=\pi_V(Z({\frak 
g}))$. In particular one has an
isomorphism
$$
{\rm End}_{U}~V\cong Z({\frak g})/Z_V({\frak g}).
$$
Note that ${\rm End}_{U}~V$ is commutative.}
\vskip 0.3cm
\noindent
{\em Proof.} Let $w\in V$ be a cyclic Whittaker vector. If $\alpha \in {\rm 
End}_{U}~V$ then $\alpha w$ is  
a Whittaker vector. But then by Theorem K there exists $u\in Z({\frak g})$ such 
that $\alpha w=uw$. 
For any $v\in U({\frak g})$ one has $\alpha vw=v\alpha w=vuw=uvw$. Thus $\alpha 
=\pi_V(u)$.

Now one can describe all irreducible Whittaker modules. A homomorphism
$$
\xi: Z({\frak g})\rightarrow {\Bbb C}
$$
is called a central character. Given a central character $\xi$ let $Z_\xi({\frak 
g})={\rm Ker}~\xi$ so that
$Z_\xi({\frak g})$ is a typical central ideal in $Z({\frak g})$.

If $V$ is any $U({\frak g})$ module one says that $V$ admits an infinitesimal 
character, and $\xi$ is its
infinitesimal character, if $\xi$ is a central character such that $uv=\xi(u)v$ 
for all $u\in Z({\frak g}),~
v\in V$. Recall that by Dixmier's theorem any irreducible $U({\frak g})$ module 
admits an 
infinitesimal character.

Given a central character $\xi$ let ${\Bbb C}_{\xi,\chi}$ be the 1-dimensional 
$Z({\frak g})\otimes U({\frak n}_+)$ module defined so that if $u\in Z({\frak 
g}),~v\in U({\frak n}_+),~
y\in {\Bbb C}_{\xi,\chi}$ then $u\otimes v y=\xi(u)\chi(v)y$. Also let
$$
Y_{\xi,\chi}=U({\frak g})\otimes_{Z({\frak g})\otimes U({\frak n}_+)}{\Bbb 
C}_{\xi,\chi} .
$$
It is clear that $Y_{\xi,\chi}$ admits an infinitesimal character and $\xi$ is 
that character.
\vskip 0.3cm
\noindent
{\bf Theorem M (\cite{K}, Theorem 3.6.1)}
{\em Let $V$ be any Whittaker module for $U({\frak g})$, the universal 
enveloping of a
simple Lie algebra $\frak g$. Then the following conditions are equivalent:

(1) $V$ is an irreducible $U({\frak g})$ module.

(2) $V$ admits an infinitesimal character.

(3) The corresponding ideal given by Theorem G is a maximal ideal.

(4) The space of Whittaker vectors in $V$ is 1--dimensional.

(5) All non--zero Whittaker vectors in $V$ are cyclic vectors.

(6) The centralizer ${\rm End}_{U}~V$ reduces to constants $\Bbb C$.

(7) $V$ is isomorphic to $Y_{\xi,\chi}$ for some central character $\xi$.}
\vskip 0.3cm
\noindent
{\em Proof.} The equivalence of (2), (3), (4) and (6) follows from Theorems K 
and L. This is also equivalent to
(5) since (5) implies that $Z({\frak g})/Z_V({\frak g})$ is a field by Theorem 
K. One gets the equivalence with 
(7) by Theorem H. It remains to relate (2)--(7) with (1). But (1) implies (2) by 
Dixmier's theorem. 
The proof of the equivalence (2)--(7) with (1) may be found in \cite{K}.


\section{Quantum deformation of the Whittaker model}\label{qWitt}

Let $\frak g$ be a complex simple Lie algebra, $U_h({\frak g})$ the standard 
quantum group associated with 
${\frak g}$.
In this section we construct a generalization of the Whittaker model $W({\frak 
b}_-)$ for $U_h({\frak g})$.

Let $U_h({\frak n}_+)$
be the subalgebra of $U_h({\frak g})$ corresponding to the nilpotent
Lie subalgebra ${\frak n}_+$. $U_h({\frak n}_+)$ is
generated by simple positive root generators of $U_h({\frak g})$
subject to the quantum Serre relations. It is easy to show that $U_h({\frak 
n}_+)$
has no non--singular characters (taking nonvanishing values
on all simple root generators). Our first main result
is a family of new realizations of the
quantum group $U_h({\frak g})$, one for each Coxeter element
in the corresponding Weyl group (see also \cite{S1}). The counterparts of 
$U({\frak n}_+)$,
which naturally arise in these new realizations of $U_h({\frak g})$,
do have non--singular characters. 

Using these new realizations we can immediately formulate a quantum group 
version
of Definition A. We also prove counterparts of Theorems A and B for $U_h({\frak 
g})$.

Finally we define quantum group generalizations of the Toda Hamiltonians. In the 
spirit of quantum harmonic
analysis these new Hamiltonians are difference operators. An alternative 
definition of these Hamiltonians has been
recently given in \cite{Et}. 


\subsection{Quantum groups}

\setcounter{equation}{0}
\setcounter{theorem}{0}

In this section we recall some basic facts about quantum groups. 
We follow the notation of \cite{ChP}. 

Let $h$ be an indeterminate, ${\Bbb C}[[h]]$ the ring of formal power series in 
$h$.
We shall consider ${\Bbb C}[[h]]$--modules equipped with the so--called 
$h$--adic 
topology. For every such module $V$ this topology is characterized by requiring 
that 
$\{ h^nV ~|~n\geq 0\}$ is a base of the neighbourhoods of $0$ in $V$, and that 
translations 
in $V$ are continuous. It is easy to see that, for modules equipped with this 
topology, every 
${\Bbb C}[[h]]$--module map is automatically continuous.

A topological Hopf algebra over ${\Bbb C}[[h]]$ is a complete ${\Bbb 
C}[[h]]$--module $A$
equipped with a structure of ${\Bbb C}[[h]]$--Hopf algebra (see \cite{ChP}, 
Definition 4.3.1),
the algebraic tensor products entering the axioms of the Hopf algebra are 
replaced by their 
completions in the $h$--adic topology.
We denote by $\mu , \imath , \Delta , \varepsilon , S$ the multiplication, the 
unit, the comultiplication,
the counit and the antipode of $A$, respectively.

The standard quantum group $U_h({\frak g})$ associated to a complex 
finite--dimensional simple Lie algebra
$\frak g$ is the algebra over ${\Bbb C}[[h]]$ topologically generated by 
elements
$H_i,~X_i^+,~X_i^-,~i=1,\ldots ,l$, and with the following defining relations:
\begin{equation}\label{qgrh}
\begin{array}{l}
[H_i,H_j]=0,~~ [H_i,X_j^\pm]=\pm a_{ij}X_j^\pm,\\
\\
X_i^+X_j^- -X_j^-X_i^+ = \delta _{i,j}\frac{K_i -K_i^{-1}}{q_i -q_i^{-1}} , \\
\\
\mbox{where }K_i=e^{d_ihH_i},~~e^h=q,~~q_i=q^{d_i}=e^{d_ih},
\end{array}
\end{equation}
and the quantum Serre relations:
$$
\begin{array}{l}
\sum_{r=0}^{1-a_{ij}}(-1)^r 
\left[ \begin{array}{c} 1-a_{ij} \\ r \end{array} \right]_{q_i} 
(X_i^\pm )^{1-a_{ij}-r}X_j^\pm(X_i^\pm)^r =0 ,~ i \neq j ,\\ \\
\mbox{ where }\\
 \\
\left[ \begin{array}{c} m \\ n \end{array} 
\right]_q=\frac{[m]_q!}{[n]_q![n-m]_q!} ,~ 
[n]_q!=[n]_q\ldots [1]_q ,~ [n]_q=\frac{q^n - q^{-n}}{q-q^{-1} }.
\end{array}
$$
$U_h({\frak g})$ is a topological Hopf algebra over ${\Bbb C}[[h]]$ with 
comultiplication 
defined by 
$$
\begin{array}{l}
\Delta_h(H_i)=H_i\otimes 1+1\otimes H_i,\\
\\
\Delta_h(X_i^+)=X_i^+\otimes K_i+1\otimes X_i^+,
\end{array}
$$
$$
\Delta_h(X_i^-)=X_i^-\otimes 1 +K_i^{-1}\otimes X_i^-,
$$
antipode defined by
$$
S_h(H_i)=-H_i,~~S_h(X_i^+)=-X_i^+K_i^{-1},~~S_h(X_i^-)=-K_iX_i^-,
$$
and counit defined by
$$
\varepsilon_h(H_i)=\varepsilon_h(X_i^\pm)=0.
$$

We shall also use the weight--type generators defined by
$$
Y_i=\sum_{j=1}^l d_i(a^{-1})_{ij}H_j,
$$ 
and the elements $L_i=e^{hY_i}$. They commute with the root vectors $X_i^\pm$ as 
follows:
\begin{equation}\label{weight-root}
L_iX_j^\pm L_i^{-1}=q_i^{\pm \delta_{ij}}X_j^\pm .
\end{equation}

The Hopf algebra $U_h({\frak g})$ is a quantization of the standard bialgebra 
structure on $\frak g$, i.e. 
$U_h({\frak g})/hU_h({\frak g})=U({\frak g}),~~ \Delta_h=\Delta~(\mbox{mod }h)$, 
where $\Delta$ is 
the standard comultiplication on $U({\frak g})$, and 
$$
\frac{\Delta_h -\Delta_h^{opp}}{h}~(\mbox{mod }h)=\delta ,
$$
where  
$\delta: {\frak g}\rightarrow {\frak g}\otimes {\frak g}$ is the standard 
cocycle on $\frak g$.
Recall that
$$
\delta (x)=({\rm ad}_x\otimes 1+1\otimes {\rm ad}_x)2r_+,~~ r_+\in {\frak 
g}\otimes {\frak g},
$$
\begin{equation}\label{rcl}
r_+=\frac 12 \sum_{i=1}^lY_i \otimes X_i + \sum_{\beta \in 
\Delta_+}(X_{\beta},X_{-\beta})^{-1} X_{\beta}\otimes X_{-\beta}.
\end{equation}
Here $X_{\pm \beta}\in {\frak g}_{\pm \beta}$ are root vectors of $\frak g$. 
The element $r_+\in {\frak g}\otimes {\frak g}$ is called a classical r--matrix.

The following proposition describes the algebraic structure of $U_h({\frak g})$. 
\begin{proposition}{\bf (\cite{ChP}, Proposition 6.5.5)}\label{algq}
Let $\frak g$ be a finite--dimensional complex simple Lie algebra, let 
$U_h({\frak h})$ be 
the subalgebra of $U_h({\frak g})$ topologically generated by the $H_i, 
i=1,\ldots l$.
Then, there is an isomorphism of algebras $\varphi :U_h({\frak g})\rightarrow 
U({\frak g})[[h]]$
over ${\Bbb C}[[h]]$ such that $\varphi =id$ (mod $h$) and $\varphi|_{U_h({\frak 
h})}=id$.
\end{proposition}

\begin{proposition}{\bf (\cite{ChP}, Proposition 6.5.7)}\label{zq}
If $\frak g$ is a finite--dimensional complex simple Lie algebra, the center 
$Z_h({\frak g})$ of 
$U_h({\frak g})$ is canonically isomorphic to $Z({\frak g})[[h]]$, where 
$Z({\frak g})$ is 
the center of $U({\frak g})$.
\end{proposition}

\begin{corollary}{\bf (\cite{ChP}, Corollary 6.5.6)}\label{rep}
If $\frak g$ be a finite--dimensional complex simple Lie algebra, then the 
assignment
$V\mapsto V[[h]]$ is a one--to--one correspondence between the 
finite--dimensional irreducible
representations of $\frak g$ and indecomposable representations of $U_h({\frak 
g})$
which are free and of finite rank as ${\Bbb C}[[h]]$--modules. Furthermore for
every such $V$ the action of the generators $H_i \in U_h({\frak g}),~~ 
i=1,\ldots l$ on 
$V[[h]]$ coincides with the action of the root generators $H_i \in {\frak h},~~ 
i=1,\ldots l$.
\end{corollary}

The representations of $U_h({\frak g})$ defined in the previous corollary are 
called
finite--dimensional representations of $U_h({\frak g})$. For every 
finite--dimensional representation
$\pi_V:{\frak g}\rightarrow {\rm End}V$ we denote the corresponding 
representation of
$U_h({\frak g})$ in the space $V[[h]]$ by the same letter.

$U_h({\frak g})$ is a quasitriangular Hopf algebra, i.e. there exists an 
invertible element
${\cal R}\in U_h({\frak g})\otimes U_h({\frak g})$, called a universal 
R--matrix, such that
\begin{equation}\label{quasitr}
\Delta^{opp}_h(a)={\cal R}\Delta_h(a){\cal R}^{-1}\mbox{ for all } a\in 
U_h({\frak g}),
\end{equation}
where $\Delta^{opp}=\sigma \Delta$, $\sigma$ is the permutation in $U_h({\frak 
g})^{\otimes 2}$,
$\sigma (x\otimes y)=y\otimes x$, and
\begin{equation}\label{rmprop}
\begin{array}{l}
(\Delta_h \otimes id){\cal R}={\cal R}_{13}{\cal R}_{23},\\
\\
(id \otimes \Delta_h){\cal R}={\cal R}_{13}{\cal R}_{12},
\end{array}
\end{equation}
where ${\cal R}_{12}={\cal R}\otimes 1,~{\cal R}_{23}=1\otimes {\cal R},
~{\cal R}_{13}=(\sigma \otimes id){\cal R}_{23}$.

From (\ref{quasitr}) and (\ref{rmprop}) it follows that ${\cal R}$ satisfies the 
quantum Yang--Baxter
equation:
\begin{equation}\label{YB}
{\cal R}_{12}{\cal R}_{13}{\cal R}_{23}={\cal R}_{23}{\cal R}_{13}{\cal R}_{12}.
\end{equation}

For every quasitriangular Hopf algebra we also have (see Proposition 4.2.7 in 
\cite{ChP}):
$$
(S\otimes id){\cal R}={\cal R}^{-1},
$$
and
\begin{equation}\label{S}
(S\otimes S){\cal R}={\cal R}.
\end{equation} 

We shall explicitly describe the element ${\cal R}$. 
First following \cite{kh-t} we recall the construction of root vectors of 
$U_h({\frak g})$. 
We shall use the so--called normal ordering in the root system  
$\Delta_+=\{\beta_1,\ldots ,\beta_N\}$ (see \cite{Z1}).
\begin{definition}\label{normord}
An ordering of the root system $\Delta_+$ is called normal if all simple roots 
are written in an arbitrary
order, and
for any theree roots $\alpha,~\beta,~\gamma$ such that
$\gamma=\alpha+\beta$ we have either $\alpha<\gamma<\beta$ or 
$\beta<\gamma<\alpha$.
\end{definition} 
To construct root vectors we shall apply the following  
inductive algorithm. Let $\alpha , \beta , \gamma \in \Delta_+$ be positive 
roots such that
$\gamma=\alpha+\beta,~\alpha<\beta$ and $[\alpha,\beta]$ is the minimal segment 
including
$\gamma$, i.e. the segment has no other roots $\alpha',\beta'$ such that 
$\gamma=\alpha'+\beta'$.
Suppose that $X_{\alpha}^\pm ,~X_{\beta}^\pm$ have
already been constructed. Then we define
\begin{equation}\label{rootvect}
\begin{array}{l}
X_{\gamma}^+=X_{\alpha}^+X_{\beta}^+ - 
q^{(\alpha,\beta)}X_{\beta}^+X_{\alpha}^+,\\
\\
X_{\gamma}^-= X_{\beta}^-X_{\alpha}^- - 
q^{-(\alpha,\beta)}X_{\alpha}^-X_{\beta}^-.
\end{array}
\end{equation} 

\begin{proposition}\label{rootprop}
For $\beta =\sum_{i=1}^lm_i\alpha_i,~m_i\in {\Bbb N}$  $X_{\beta}^\pm $ is a 
polynomial in 
the noncommutative variables $X_i^\pm$ homogeneous in each $X_i^\pm$ of degree 
$m_i$.
\end{proposition}

The root vectors $X_{\beta}$ satisfy the following relations:
$$
[X_\alpha^+,X_{\alpha}^-]=a(\alpha)\frac{e^{h\alpha^\vee}-e^{-h\alpha^\vee}}{ 
q-q^{-1}}.
$$
where $a(\alpha)\in {\Bbb C}[[h]]$.
They commute with elements of the subalgebra $U_h({\frak h})$ as follows:

\begin{equation}\label{roots-cart}
[H_i,X_{\beta}^\pm]=\pm \beta(H_i)X_{\beta}^\pm,~i=1,\ldots ,l.
\end{equation}

Note that by construction
$$
\begin{array}{l}
X_\beta^+~(\mbox{mod }h)=X_\beta \in {\frak g}_\beta,\\
\\
X_\beta^-~(\mbox{mod }h)=X_{-\beta} \in {\frak g}_{-\beta}
\end{array}
$$
are root vectors of $\frak g$. This implies that $a(\alpha)~(\mbox{mod 
}h)=(X_\alpha,X_{-\alpha})$.

Let $U_h({\frak n}_+),U_h({\frak n}_-)$ be the ${\Bbb C}[[h]]$--subalgebras of 
$U_h({\frak g})$ topologically 
generated by the 
$X_i^+$ and by the $X_i^-$, respectively.

Now using the root vectors $X_{\beta}^\pm$ we can construct a topological basis 
of 
$U_h({\frak g})$.
Define for ${\bf r}=(r_1,\ldots ,r_N)\in {\Bbb N}^N$,
$$
(X^+)^{\bf r}=(X_{\beta_1}^+)^{r_1}\ldots (X_{\beta_N}^+)^{r_N},
$$
$$
(X^-)^{\bf r}=(X_{\beta_1}^-)^{r_1}\ldots (X_{\beta_N}^-)^{r_N},
$$
and for ${\bf s}=(s_1,\ldots s_l)\in {\Bbb N}^{~l}$,
$$
H^{\bf s}=H_1^{s_1}\ldots H_l^{s_l}.
$$
\begin{proposition}{\bf (\cite{kh-t}, Proposition 3.3)}\label{PBW}
The elements $(X^+)^{\bf r}$, $(X^-)^{\bf t}$ and $H^{\bf s}$, for ${\bf 
r},~{\bf t}\in {\Bbb N}^N$, 
${\bf s}\in {\Bbb N}^l$, form topological bases of $U_h({\frak n}_+),U_h({\frak 
n}_-)$ and $U_h({\frak h})$,
respectively, and the products $(X^+)^{\bf r}H^{\bf s}(X^-)^{\bf t}$ form a 
topological basis of 
$U_h({\frak g})$. In particular, multiplication defines an isomorphism of ${\Bbb 
C}[[h]]$ modules:
$$
U_h({\frak n}_-)\otimes U_h({\frak h}) \otimes U_h({\frak n}_+)\rightarrow 
U_h({\frak g}).
$$
\end{proposition}

An explicit expression for ${\cal R}$ may be written by making use of the 
q-exponential
$$
exp_q(x)=\sum_{k=0}^\infty \frac{x^k}{(k)_q!},
$$ 
where
$$
(k)_q!=(1)_q\ldots (k)_q,~~(n)_q=\frac{q^n -1}{q-1}.
$$

Now the element ${\cal R}$ may be written as (see Theorem 8.1 in \cite{kh-t}):
\begin{equation}\label{univr}
{\cal R}=exp\left[ h\sum_{i=1}^l(Y_i\otimes H_i)\right]\prod_{\beta}
exp_{q_{\beta}^{-1}}[(q-q^{-1})a(\beta)^{-1}X_{\beta}^+\otimes X_{\beta}^-],
\end{equation}
where $q_\beta =q^{(\beta,\beta)}$;
the product is over all the positive roots of $\frak g$, and the order of the 
terms is such that 
the $\alpha$--term appears to the left of the $\beta$--term if $\alpha <\beta$ 
with respect to the normal
ordering of $\Delta_+$.

\begin{remark}
The r--matrix $r_+=\frac 12 h^{-1}({\cal R}-1\otimes 1)~~(\mbox{mod }h)$, which 
is the classical limit of ${\cal R}$,
coincides with the classical r--matrix (\ref{rcl}).
\end{remark}


\subsection{Non--singular characters and quantum groups}

\setcounter{equation}{0}
\setcounter{theorem}{0}

In this section following \cite{S1} we recall the construction of quantum 
counterparts of the principal nilpotent 
Lie subalgebras of complex simple Lie algebras and of their non--singular 
characters. Subalgebras of $U_h({\frak g})$ which resemble the subalgebra 
$U({\frak n}_+) 
\subset U({\frak g})$ and possess non--singular
characters naturally appear in the Coxeter realizations of $U_h({\frak g})$
defined in \cite{S1} as follows.

Denote by $S_l$ the symmetric group of $l$ elements.
To any element $\pi \in S_l$ we associate a Coxeter element $s_{\pi}$ by the 
formula
$s_\pi =s_{\pi (1)}\ldots s_{\pi (l)}$.
Let 
$U_h^{s_\pi}({\frak g})$ be the associative algebra over ${\Bbb C}[[h]]$ with 
generators 
$e_i , f_i , H_i,~i=1, \ldots l$ subject to the relations:
\begin{equation}\label{sqgr}
\begin{array}{l}
[H_i,H_j]=0,~~ [H_i,e_j]=a_{ij}e_j, ~~ [H_i,f_j]=-a_{ij}f_j,\\
\\
e_i f_j -q^{ c^\pi _{ij}} f_j e_i = \delta _{i,j}\frac{K_i -K_i^{-1}}{q_i 
-q_i^{-1}} , \\
 \\
K_i=e^{d_ihH_i}, \\
 \\
\sum_{r=0}^{1-a_{ij}}(-1)^r q^{r c_{ij}^\pi}
\left[ \begin{array}{c} 1-a_{ij} \\ r \end{array} \right]_{q_i} 
(e_i )^{1-a_{ij}-r}e_j (e_i)^r =0 ,~ i \neq j , \\
 \\
\sum_{r=0}^{1-a_{ij}}(-1)^r q^{r c_{ij}^\pi}
\left[ \begin{array}{c} 1-a_{ij} \\ r \end{array} \right]_{q_i} 
(f_i )^{1-a_{ij}-r}f_j (f_i)^r =0 ,~ i \neq j ,
\end{array}
\end{equation}
where $c_{ij}^{\pi}=\left( \frac{1+s_\pi}{1-s_\pi }\alpha_i , \alpha_j \right)$ 
are matrix elements of the Cayley transform of $s_\pi$ 
in the basis of simple roots.

\begin{theorem}{\bf (\cite{S1}, Theorem 4)}\label{newreal}
For every solution $n_{ij}\in {\Bbb C},~i,j=1,\ldots ,l$ of equations 
\begin{equation}\label{eqpi}
d_jn_{ij}-d_in_{ji}=c^\pi_{ij}
\end{equation}
there exists an algebra
isomorphism $\psi_{\{ n\}} : U_h^{s_\pi}({\frak g}) \rightarrow 
U_h({\frak g})$ defined  by the formulas:
$$
\begin{array}{l}
\psi_{\{ n\}}(e_i)=X_i^+ \prod_{p=1}^lL_p^{n_{ip}},\\
 \\
\psi_{\{ n\}}(f_i)=\prod_{p=1}^lL_p^{-n_{ip}}X_i^- , \\
 \\
\psi_{\{ n\}}(H_i)=H_i .
\end{array}
$$
\end{theorem}

We call the algebra $U_h^{s_\pi}({\frak g})$ the Coxeter realization of the 
quantum group $U_h({\frak g})$ corresponding to the Coxeter
element $s_\pi$.

Let $U_h^{s_\pi}({\frak n}_+) $ be the subalgebra in $U_h^{s_\pi}({\frak g})$ 
generated by 
$e_i ,i=1, \ldots l$.

\begin{proposition}{\bf (\cite{S1}, Proposition 2)}\label{charf}
The map $\chi_h^{s_\pi}:U_h^{s_\pi}({\frak n}_+) \rightarrow {\Bbb C}[[h]]$ 
defined on generators by 
$\chi_h^{s_\pi}(e_i)=c_i,~c_i\in {\Bbb C}[[h]],~c_i\neq 0$ is a character of the 
algebra $U_h^{s_\pi}({\frak n}_+) $.
\end{proposition}

The proof of this proposition given in \cite{S1} is based on the following 
Lemma.

\begin{lemma}{\bf (\cite{S1}, Lemma 3)}\label{tmatrel}
The matrix elements of  $\frac{1+s_\pi}{1-s_\pi }$ are of the form :
\begin{equation}\label{matrel}
c_{ij}^{\pi}=\left( \frac{1+s_\pi}{1-s_\pi }\alpha_i , \alpha_j \right)=
\varepsilon_{ij}^\pi b_{ij},
\end{equation}
where
$$
\varepsilon_{ij}^\pi =\left\{ \begin{array}{ll}
-1 & \pi^{-1}(i) <\pi^{-1}(j) \\
0 & i=j \\
1 & \pi^{-1}(i) >\pi^{-1}(j) 
\end{array}
\right  .
$$
\end{lemma}

Now we shall study the algebraic structure of $U_h^{s_\pi}({\frak g})$.
Denote by $U_h^{s_\pi}({\frak n}_-) $ the subalgebra in $U_h^{s_\pi}({\frak g})$ 
generated by
$f_i ,i=1, \ldots l$. From defining relations (\ref{sqgr}) it follows that the 
map 
$\overline \chi_h^{s_\pi}:U_h^{s_\pi}({\frak n}_-) \rightarrow {\Bbb C}[[h]]$ 
defined on generators by 
$\overline \chi_h^{s_\pi}(f_i)=c_i, c_i\in {\Bbb C}[[h]], c_i\neq 0$ is a 
character of the algebra $U_h^{s_\pi}({\frak n}_-)$.

Let $U_h^{s_\pi}({\frak h})$ be the subalgebra in $U_h^{s_\pi}({\frak g})$ 
generated by $H_i,~i=1,\ldots ,l$.
Define $U_h^{s_\pi}({\frak b}_\pm)=U_h^{s_\pi}({\frak n}_\pm)U_h^{s_\pi}({\frak 
h})$.

We shall construct a Poincar\'{e}--Birkhoff-Witt basis for $U_h^{s_\pi}({\frak 
g})$.
It is convenient to introduce an operator $K\in {\rm End}~{\frak h}$ such that
\begin{equation}\label{Kdef}
KH_i=\sum_{j=1}^l\frac{n_{ij}}{d_i}Y_j.
\end{equation}
In particular, we have 
$$
\frac{n_{ji}}{d_j}=(KH_j,H_i).
$$

Equation (\ref{eqpi}) is equivalent to the following equation for the operator 
$K$:
$$
K-K^* = \frac{1+s_\pi}{1-s_\pi}.
$$

\begin{proposition}\label{rootss}
(i)For any solution of equation (\ref{eqpi}) and any normal ordering of the root 
system $\Delta_+$
the elements $e_{\beta}=\psi_{\{ n\}}^{-1}(X_{\beta}^+e^{hK\beta^\vee})$ and 
$f_{\beta}=\psi_{\{ n\}}^{-1}(e^{-hK\beta^\vee}X_{\beta}^-),~\beta \in \Delta_+$
lie in the subalgebras $U_h^{s_\pi}({\frak n}_+)$ and $U_h^{s_\pi}({\frak 
n}_-)$, respectively.

(ii)Moreover, the elements 
$e^{\bf r}=e_{\beta_1}^{r_1}\ldots e_{\beta_N}^{r_N},~~f^{\bf 
t}=e_{\beta_1}^{t_1}\ldots e_{\beta_N}^{t_N}$
and $H^{\bf s}=H_1^{s_1}\ldots H_l^{s_l}$
for ${\bf r},~{\bf t},~{\bf s}\in {\Bbb N}^N$, form 
topological bases of $U_h^{s_\pi}({\frak n}_+),~U_h^{s_\pi}({\frak n}_-)$ and 
$U_h^{s_\pi}({\frak h})$, 
and the products $f^{\bf t}H^{\bf s}e^{\bf r}$ form a topological basis of 
$U_h^{s_\pi}({\frak g})$. In particular, multiplication defines an isomorphism 
of ${\Bbb C}[[h]]$ modules$$
U_h^{s_\pi}({\frak n}_-)\otimes U_h^{s_\pi}({\frak h})\otimes U_h^{s_\pi}({\frak 
n}_+)\rightarrow U_h^{s_\pi}({\frak g}).
$$
\end{proposition}
{\em Proof.} Let $\beta=\sum_{i=1}^l m_i\alpha_i \in \Delta_+$ be a positive 
root, 
$X_{\beta}^+\in U_h({\frak g})$ the corresponding root vector. Then 
$\beta^\vee=\sum_{i=1}^l m_id_iH_i$, and so
$K\beta^\vee=\sum_{i,j=1}^l m_in_{ij}Y_j$. Now the proof of the first statement 
follows immediately from 
Proposition \ref{rootprop}, commutation relations (\ref{weight-root}) and the 
definition of the isomorphism
$\psi_{\{ n\}}$. The second assertion is a consequence of Proposition \ref{PBW}.

Now we would like to choose a normal ordering of the root system $\Delta_+$ in 
such a way that 
$\chi_h^{s_\pi}(e_{\beta})=0$ and $\overline \chi_h^{s_\pi}(f_{\beta})=0$ if 
$\beta$ is not a simple root.  
\begin{proposition}\label{rootsh}
Choose a normal ordering of the root system $\Delta_+$ such that the simple 
roots are written
in the following order: $\alpha_{\pi (1)},\ldots ,\alpha_{\pi (l)}$.  
Then $\chi_h^{s_\pi}(e_{\beta})=0$ and $\overline \chi_h^{s_\pi}(f_{\beta})=0$ 
if $\beta$ is not a simple root.
\end{proposition} 
{\em Proof.} We shall consider the case of positive root generators.
The proof for 
negative root generators is similar to that for the positive ones.
 
The root vectors $X_{\beta}^+$ are defined in terms of iterated q-commutators
(see (\ref{rootvect})). Therefore it suffices to verify that for $i<j$
$$
\begin{array}{l}
\chi_h^{s_\pi}(e_{\alpha_{\pi(i)}+\alpha_{\pi(j)}})=\\
\\
\chi_h^{s_\pi}(\psi_{\{ n\}}^{-1}( (X_{\pi(i)}^+X_{\pi(j)}^+ - 
q^{(\alpha_{\pi(i)},\alpha_{\pi(j)})}X_{\pi(j)}^+X_{\pi(i)}^+)
e^{hK(d_{\pi(i)}H_{\pi(i)}+d_{\pi(j)}H_{\pi(j)})}))=0.
\end{array}
$$

From (\ref{Kdef}) and commutation relations (\ref{weight-root}) we obtain that
\begin{equation}\label{bebe}
\begin{array}{l}
\psi_{\{ n\}}^{-1}((X_{\pi(i)}^+X_{\pi(j)}^+ - 
q^{(\alpha_{\pi(i)},\alpha_{\pi(j)})}X_{\pi(j)}^+X_{\pi(i)}^+)
e^{hK(d_{\pi(i)}H_{\pi(i)}+d_{\pi(j)}H_{\pi(j)})})= \\
\\
q^{-d_{\pi(j)}n_{\pi(i)\pi(j)}}(e_{\pi(i)}e_{\pi(j)} - 
q^{b_{\pi(i)\pi(j)}+d_{\pi(j)}n_{\pi(i)\pi(j)}-d_{\pi(i)}n_{\pi(j)\pi(i)}}e_{\pi
(j)}e_{\pi(i)})
\end{array}
\end{equation}

Now using equation (\ref{eqpi}) and Lemma \ref{tmatrel} the combination 
$b_{\pi(i)\pi(j)}+d_{\pi(j)}n_{\pi(i)\pi(j)}-d_{\pi(i)}n_{\pi(j)\pi(i)}$ may be 
represented as
$b_{\pi(i)\pi(j)}+\varepsilon_{\pi(i)\pi(j)}^\pi b_{\pi(i)\pi(j)}$. But 
$\varepsilon_{\pi(i)\pi(j)}^\pi =-1$ for $i<j$ and therefore the r.h.s. of 
(\ref{bebe}) takes the form
$$
q^{-d_{\pi(j)}n_{\pi(i)\pi(j)}}[e_{\pi(i)},e_{\pi(j)}].
$$
Clearly,
$$
\chi_h^{s_\pi}(e_{\alpha_{\pi(i)}+\alpha_{\pi(j)}})=
q^{-d_{\pi(j)}n_{\pi(i)\pi(j)}}\chi_h^{s_\pi}([e_{\pi(i)},e_{\pi(j)}])=0.
$$


\subsection{Quantum deformation of the Whittaker model}\label{whitth}

\setcounter{equation}{0}
\setcounter{theorem}{0}

In this section we define a quantum deformation of the Whittaker model $W({\frak 
b}_-)$.
Our construction is similar the one described in Section \ref{whitt}, the 
quantum group $U_h^{s_\pi}({\frak g})$,
the subalgebra $U_h^{s_\pi}({\frak n}_+)$ and characters 
${\chi_h^{s_\pi}}:U_h^{s_\pi}({\frak n}_+) \rightarrow {\Bbb C}[[h]]$ serve as 
natural
counterparts of the universal enveloping algebra $U({\frak g})$,
of the subalgebra $U({\frak n}_+)$ and of non--singular characters
$\chi:U({\frak n}_+) \rightarrow {\Bbb C}$.

Let ${U_h^{s_\pi}({\frak n}_+)}_{\chi_h^{s_\pi}}$ be the kernel of the character
$\chi_h^{s_\pi}:U_h^{s_\pi}({\frak n}_+) \rightarrow {\Bbb C}[[h]]$ so that
one has a direct sum
$$
U_h^{s_\pi}({\frak n}_+)={\Bbb C}[[h]]\oplus {U_h^{s_\pi}({\frak 
n}_+)}_{\chi_h^{s_\pi}}.
$$

From Proposition \ref{rootss} we have a linear 
isomorphism $U_h^{s_\pi}({\frak g})=U_h^{s_\pi}({\frak b}_-)\otimes 
U_h^{s_\pi}({\frak n}_+)$ and hence 
the direct sum
\begin{equation}\label{maindecq}
U_h^{s_\pi}({\frak g})=U_h^{s_\pi}({\frak b}_-) \oplus I_{\chi_h^{s_\pi}},
\end{equation}
where $I_{\chi_h^{s_\pi}}=U_h^{s_\pi}({\frak g}){U_h^{s_\pi}({\frak 
n}_+)}_{\chi_h^{s_\pi}}$ is the left--sided ideal generated by 
${U_h^{s_\pi}({\frak n}_+)}_{\chi_h^{s_\pi}}$.

For any $u\in U_h^{s_\pi}({\frak g})$ let $u^{\chi_h^{s_\pi}}\in 
U_h^{s_\pi}({\frak b}_-)$ be its component in 
$U_h^{s_\pi}({\frak b}_-)$ relative to the decomposition (\ref{maindecq}). 
Denote by $\rho_{\chi_h^{s_\pi}}$
the linear map 
$$
\rho_{\chi_h^{s_\pi}} : U_h^{s_\pi}({\frak g}) \rightarrow U_h^{s_\pi}({\frak 
b}_-)
$$
given by $\rho_{\chi_h^{s_\pi}} (u)=u^{\chi_h^{s_\pi}}$.

Denote by $Z_h^{s_\pi}({\frak g})$ the center of $U_h^{s_\pi}({\frak g})$.
From Proposition \ref{zq} and Theorem \ref{newreal} we obtain that 
$Z_h^{s_\pi}({\frak g})\cong Z({\frak g})[[h]]$. In particular, 
$Z_h^{s_\pi}({\frak g})$
is freely generated as a commutative topological algebra over ${\Bbb C}[[h]]$ by 
$l$ elements $I_1,\ldots , I_l$.

Let $W_h({\frak b}_-)=\rho_{\chi_h^{s_\pi}} (Z_h^{s_\pi}({\frak g}))$. 
\vskip 0.3cm
\noindent
{\bf Theorem $\bf A_h$ }
{\em The map
\begin{equation}\label{mapq}
\rho_{\chi_h^{s_\pi}} : Z_h^{s_\pi}({\frak g}) \rightarrow W_h({\frak b}_-)
\end{equation}
is an isomorphism of algebras. In particular, $W_h({\frak b}_-)$  
is freely generated as a commutative topological algebra over ${\Bbb C}[[h]]$ by 
$l$ elements 
$I_i^{\chi_h^{s_\pi}}=\rho_{\chi_h^{s_\pi}}(I_i),~~i=1,\ldots ,l$.}
\vskip 0.3cm
\noindent
{\em Proof} is similar to that of Theorem A in the classical case.
\begin{remark}\label{homrhoh}
Similarly to Remark \ref{homrho} we have 
$$
\rho_\chi(uv)=\rho_\chi(u)\rho_\chi(v)
$$
for any $u\in U({\frak g}),~~v\in Z({\frak g})$.
\end{remark} 

\noindent
{\bf Definition $\bf A_h$ }
{\em The algebra $W_h({\frak b}_-)$ is called the Whittaker model of 
$Z_h^{s_\pi}({\frak g})$.}  
\vskip 0.3cm

Next we equip $U_h^{s_\pi}({\frak b}_-)$ with a structure of a left 
$U_h^{s_\pi}({\frak n}_+)$ module in such a
way that $W_h({\frak b}_-)$ is identified with the space of invariants with 
respect to this action.
Following Lemma A in the classical case we define this action by
\begin{equation}\label{mainactq}
x\cdot v =[x,v]^{\chi_h^{s_\pi}},
\end{equation}
where $v\in U_h^{s_\pi}({\frak b}_-)$ and $x\in U_h^{s_\pi}({\frak n}_+)$.

Consider the space $U_h^{s_\pi}({\frak b}_-)^{U_h^{s_\pi}({\frak n}_+)}$ of 
$U_h^{s_\pi}({\frak n}_+)$ invariants 
in $U_h^{s_\pi}({\frak b}_-)$ with respect to this
action. Clearly, 
$W_h({\frak b}_-)\subseteq U_h^{s_\pi}({\frak b}_-)^{U_h^{s_\pi}({\frak n}_+)}$.

\noindent
{\bf Theorem $\bf B_h$ }
{\em Suppose that $\chi_h^{s_\pi}(e_i)\neq 0~(\mbox{mod }h)$ for $i=1,\ldots l$. 
Then the space of $U_h^{s_\pi}({\frak n}_+)$ invariants 
in $U_h^{s_\pi}({\frak b}_-)$ with respect to the
action (\ref{mainactq}) is isomorphic to $W_h({\frak b}_-)$, i.e.}
\begin{equation}\label{invq}
U_h^{s_\pi}({\frak b}_-)^{U_h^{s_\pi}({\frak n}_+)}\cong W_h({\frak b}_-).
\end{equation}

The proof of this theorem, as well as proofs of many statements in this paper, 
is based on the following Lemma.

\begin{lemma}\label{mainl}
Let $V$ be a complete ${\Bbb C}[[h]]$ module, $A, B\subset V$ two closed 
subspaces. Denote by 
$p:V\rightarrow V/hV$ the canonical projection. Suppose that
$B\subseteq A$, $p(A)=p(B)$ and for any $a\in A,~b\in B$ such that $a-b=hc,~c\in 
V$
we have $c\in A$. Then $A=B$.
\end{lemma} 
\noindent
{\em Proof.}
Let $a\in A$. Since $p(A)=p(B)$ one can find an element $b_0\in B$ such
that $a-b_0=ha_1,~a_1\in A$. Applying the same procedure
to $a_1$ one can find elements $b_1\in B,
~a_2\in A$ such that $a_1-b_1=ha_2$, i.e. 
$a-b_0-hb_1=0~(\mbox{mod }h^2)$. 
We can continue this process. Finally we obtain an infinite sequence of elements 
$b_i\in B$ such that $a-\sum_{i=0}^p h^pb_p=0~(\mbox{mod }h^{p+1})$. Since the 
subspace $A$ 
is closed in the $h$--adic topology
the series 
$\sum_{i=0}^\infty h^pb_p\in B$ 
converges to $a$. Therefore $a\in B$. This completes the proof.

\vskip 0.3cm
\noindent
{\em Proof of Theorem $B_h$.}
Let $p: U_h^{s_\pi}({\frak g})\rightarrow U_h^{s_\pi}({\frak 
g})/hU_h^{s_\pi}({\frak g})=U({\frak g})$ 
be the canonical projection. Note that $p(U_h^{s_\pi}({\frak n}_+))=U({\frak 
n}_+),~
p(U_h^{s_\pi}({\frak b}_-))=U({\frak b}_-)$ and for every $x\in 
U_h^{s_\pi}({\frak n}_+)$  
$\chi_h^{s_\pi}(x)~(\mbox{mod }h)=\chi(p(x))$ for some non--singular character 
$\chi: U({\frak n}_+)\rightarrow {\Bbb C}$. 
Therefore $p(\rho_{\chi_h^{s_\pi}}(x))=\rho_\chi(p(x))$ for every $x\in 
U_h^{s_\pi}({\frak g})$, 
and hence by Theorem ${\rm A}_q$ $p(W_h({\frak b}_-))=W({\frak b}_-)$. Using 
Lemma A and 
the definition of action (\ref{mainactq}) we also obtain that 
$p(U_h^{s_\pi}({\frak b}_-)^{U_h^{s_\pi}({\frak n}_+)})=
U({\frak b}_-)^{N_+} = W({\frak b}_-)$.

Now Theorem ${\rm B_h}$ follows immediately from Lemma \ref{mainl} applied to 
$V=U_h^{s_\pi}({\frak g}),~~
A=U_h^{s_\pi}({\frak b}_-)^{U_h^{s_\pi}({\frak n}_+)},~~B=W_h({\frak b}_-)$.


\subsection{Quantum deformations of Whittaker modules}\label{Whittmodh}

\setcounter{equation}{0}
\setcounter{theorem}{0}

In this section we define quantum deformations of Whittaker modules. 
The construction of these modules
is similar to that for Lie algebras (see Section \ref{Whittmod}).

Fix a Coxeter element $s_\pi\in W$.
Let $V_h$ be a $U_h^{s_\pi}({\frak g})$ module, which is also free as a ${\Bbb 
C}[[h]]$ module. 
The action is denoted by $uv$ for $u\in U_h^{s_\pi}({\frak g}),~~v\in V_h$. Let 
${\chi_h^{s_\pi}}:U_h^{s_\pi}({\frak n}_+) \rightarrow {\Bbb C}[[h]]$ be a 
non--singular character of $U_h^{s_\pi}({\frak n}_+)$ 
(see Proposition \ref{charf}). We shall assume that ${\chi_h^{s_\pi}}(e_i)\neq 
0\mbox{ mod h},~i=1,\ldots l$.
We call a vector $w_h\in V_h$ a Whittaker vector (with respect to 
${\chi_h^{s_\pi}}$) if
$$
xw_h={\chi_h^{s_\pi}} (x)w_h
$$
for all $x\in U_h^{s_\pi}({\frak n}_+)$. A Whittaker vector $w_h$ is called a 
cyclic Whittaker vector (for $V_h$) if 
$U_h^{s_\pi}({\frak g})w_h=V_h$. A $U_h^{s_\pi}({\frak g})$ module $V_h$ is 
called a Whittaker module if it contains a cyclic 
Whittaker vector. 
\begin{remark}\label{whitto}
Note that in this case $V=V_h/hV_h$ is naturally a Whittaker module for 
$U({\frak g})$, 
$w=w_h\mbox{ mod h}\in V$ being a cyclic Whittaker vector for $V$ with respect 
to the non--singular character $\chi$ of 
$U({\frak n}_+)$ defined by $\chi(e_i\mbox{ mod h})={\chi_h^{s_\pi}}(e_i)\mbox{ 
mod h}$.
\end{remark}

If $V_h$ is any $U_h^{s_\pi}({\frak g})$ module we let $U_{h,V}^{s_\pi}({\frak 
g})$ be the annihilator of $V_h$. Then $U_{h,V}^{s_\pi}({\frak g})$
defines a central ideal $Z_{h,V}^{s_\pi}({\frak g})$ by putting
\begin{equation}\label{zvh}
Z_{h,V}^{s_\pi}({\frak g})=Z_h^{s_\pi}({\frak g})\cap U_{h,V}^{s_\pi}({\frak 
g}).
\end{equation}

Now assume that $V_h$ is a Whittaker module for $U_h^{s_\pi}({\frak g})$ and 
$w_h\in V_h$ is a cyclic Whittaker vector.
Let $U_{h,w}^{s_\pi}({\frak g})\subseteq U_h^{s_\pi}({\frak g})$ be the 
annihilator of $w_h$. Thus $U_{h,V}^{s_\pi}({\frak g}) \subseteq
U_{h,w}^{s_\pi}({\frak g})$, where $U_{h,w}^{s_\pi}({\frak g})$ is a left ideal 
and $U_{h,V}^{s_\pi}({\frak g})$ is a two--sided ideal in 
$U_h^{s_\pi}({\frak g})$. One has $V_h=U_h^{s_\pi}({\frak 
g})/U_{h,w}^{s_\pi}({\frak g})$ as $U_h^{s_\pi}({\frak g})$ modules so that 
$V_h$ is determined
up to equivalence by $U_{h,w}^{s_\pi}({\frak g})$. Clearly 
$I_{\chi_h^{s_\pi}}=U_h^{s_\pi}({\frak g}){U_h^{s_\pi}({\frak 
n}_+)}_{\chi_h^{s_\pi}}\subseteq 
U_{h,w}^{s_\pi}({\frak g})$ and $U_h^{s_\pi}({\frak g})Z_{h,V}^{s_\pi}({\frak 
g})\subseteq U_{h,w}^{s_\pi}({\frak g})$.

The following theorem, similar to Theorem F for Lie algebras, says that up to 
equivalence $V_h$ is determined by the central ideal $Z_{h,V}^{s_\pi}({\frak 
g})$.
\vskip 0.3cm
\noindent
{\bf Theorem $\bf F_h$ }
{\em Let $V_h$ be any $U_h^{s_\pi}({\frak g})$ module which admits a cyclic 
Whittaker vector $w_h$ and let
$U_{h,w}^{s_\pi}({\frak g})$ be the annihilator of $w_h$. Then}
\begin{equation}\label{anndech}
U_{h,w}^{s_\pi}({\frak g})=U_h^{s_\pi}({\frak g})Z_{h,V}^{s_\pi}({\frak 
g})+I_{\chi_h^{s_\pi}} .
\end{equation}
\vskip 0.3cm
\noindent
{\em Proof.} Denote by 
$p: U_h^{s_\pi}({\frak g})\rightarrow U_h^{s_\pi}({\frak g})/hU_h^{s_\pi}({\frak 
g})\cong U({\frak g})$ the canonical projection.
Then we have $p(U_{h,w}^{s_\pi}({\frak g}))=U_{w}({\frak g})$, where 
$U_{w}({\frak g})$ is the annihilator of the 
Whittaker vector $w=w_h\mbox{ mod h}$ in the Whittaker module $V=V_h/hV_h$ for 
$U({\frak g})$ (see Remark \ref{whitto}).
On the other hand from Theorem F we have $U_w({\frak g})=U({\frak g})Z_V({\frak 
g})+I_\chi$. But clearly
$p(U_h^{s_\pi}({\frak g})Z_{h,V}^{s_\pi}({\frak g})+I_{\chi_h^{s_\pi}})=U({\frak 
g})Z_V({\frak g})+I_\chi$.
Now the result of the theorem follows immediately from Lemma \ref{mainl} applied 
to $A=U_{h,w}^{s_\pi}({\frak g})$ and 
$B=U_h^{s_\pi}({\frak g})Z_{h,V}^{s_\pi}({\frak g})+I_{\chi_h^{s_\pi}}$.
\vskip 0.3cm
We also have the following quantum counterpart of Lemma B. We use the notation 
of Section \ref{whitth}.
If $X\subseteq U_h^{s_\pi}({\frak g})$ let $X^{\chi_h^{s_\pi}} 
=\rho_{\chi_h^{s_\pi}}(X)$. Note that $U_{h,w}^{s_\pi}({\frak g})$ is stable 
under the map
$u\mapsto \rho_{\chi_h^{s_\pi}}(u)$. We recall also that by Theorem $\rm A_h$ 
$\rho_{\chi_h^{s_\pi}}$ induces an algebra isomorphism
$Z_h^{s_\pi}({\frak g}) \rightarrow W_{h}({\frak b}_-)$, where $W_{h}({\frak 
b}_-)=Z_h^{s_\pi}({\frak g})^{\chi_h^{s_\pi}}$. Thus if $Z_{h,*}$ is any ideal 
in
$Z_h^{s_\pi}({\frak g})$ then $W_{h,*}({\frak b}_-)=Z_{h,*}^{\chi_h^{s_\pi}}$ is 
an isomorphic ideal in $W_{h}({\frak b}_-)$.
But $(U_h^{s_\pi}({\frak g})Z_{h,*})^{\chi_h^{s_\pi}}=U_{h}^{s_\pi}({\frak 
b}_-)W_{h,*}({\frak b}_-)$ by Remark \ref{homrhoh}. Thus by (\ref{maindecq})
one has the direct sum
\begin{equation}\label{anndec*h}
U_h^{s_\pi}({\frak g})Z_{h,*}+I_{\chi_h^{s_\pi}}=U_{h}^{s_\pi}({\frak 
b}_-)W_{h,*}({\frak b}_-)\oplus I_{\chi_h^{s_\pi}}.
\end{equation}
\vskip 0.3cm
\noindent
{\bf Lemma $\bf B_h$}
{\em Let $X=\{ v \in U_{h}^{s_\pi}({\frak b}_-) | (x\cdot v)w_h=0 \mbox{ for all 
}x \in U_h^{s_\pi}({\frak n}_+)\}$, where
$x\cdot v$ is given by (\ref{mainactq}). Then
\begin{equation}
X=U_{h}^{s_\pi}({\frak b}_-)W_{h,V}({\frak b}_-)+W_{h}({\frak b}_-),
\end{equation}
where $W_{h,V}({\frak b}_-)=Z_{h,V}^{s_\pi}({\frak g})^{\chi_h^{s_\pi}}$. 
Furthermore if we denote 
$U_{h,w}^{s_\pi}({\frak b}_-)=U_{h,w}^{s_\pi}({\frak g})\cap 
U_{h}^{s_\pi}({\frak b}_-)$ then
\begin{equation}\label{bannh}
U_{h,w}^{s_\pi}({\frak b}_-)=U_{h}^{s_\pi}({\frak b}_-)W_{h,V}({\frak b}_-).
\end{equation}
}
\vskip 0.3cm
\noindent
{\em Proof} is similar to that of Theorem $\rm F_h$: one should apply Lemma 
\ref{mainl} to $A=X$ and 
$B=U_{h}^{s_\pi}({\frak b}_-)W_{h,V}({\frak b}_-)+W_{h}({\frak b}_-)$.
\vskip 0.3cm

Now one can determine, up to equivalence, the set of all Whittaker modules for 
$U_h^{s_\pi}({\frak g})$. 
They are naturally
parametrized by the set of all ideals in the center $Z_h^{s_\pi}({\frak g})$.
\begin{remark}
The proofs of Theorems $\rm G_h$, $\rm H_h$, $\rm K_h$ and $\rm L_h$ below  are 
based on Theorems $\rm A_h$, $\rm F_h$,
Lemma $\rm B_h$ and completely similar to the proofs of Theorems $\rm G$, $\rm 
H$, $\rm K$ and $\rm L$ in the classical case
(see Section \ref{Whittmod}). We do not reproduce these proofs in this section.
\end{remark}
\noindent
{\bf Theorem $\bf G_h$} 
{\em Let $V_h$ be any Whittaker module for $U_h^{s_\pi}({\frak g})$. Let 
$U_{h,V}^{s_\pi}({\frak g})$ be the annihilator of $V_h$ and let 
$Z_h^{s_\pi}({\frak g})$
be the center of $U_h^{s_\pi}({\frak g})$. Then the correspondence
\begin{equation}\label{idcorrh}
V_h\mapsto Z_{h,V}^{s_\pi}({\frak g}),
\end{equation}
where $Z_{h,V}^{s_\pi}({\frak g})=U_{h,V}^{s_\pi}({\frak g})\cap 
Z_h^{s_\pi}({\frak g})$, sets up a bijection between the set of all 
equivalence classes of Whittaker modules and the set of all ideals in 
$Z_h^{s_\pi}({\frak g})$.}
\vskip 0.3cm
\noindent

Now consider the subalgebra $Z_h^{s_\pi}({\frak g})U_h^{s_\pi}({\frak n}_+)$ in 
$U_h^{s_\pi}({\frak g})$. Using the arguements applied in the proof of 
Theorem $\rm A_h$ it is aesy to show that
$Z_h^{s_\pi}({\frak g})U_h^{s_\pi}({\frak n}_+)\cong Z_h^{s_\pi}({\frak 
g})\otimes U_h^{s_\pi}({\frak n}_+)$. Now let $Z_{h,*}$ be any ideal in
$Z_h^{s_\pi}({\frak g})$ and regard $Z_h^{s_\pi}({\frak g})/Z_{h,*}$ as a 
$Z_h^{s_\pi}({\frak g})$ module. Equip $Z_h^{s_\pi}({\frak g})/Z_{h,*}$ with a 
structure of $Z_h^{s_\pi}({\frak g})\otimes U_h^{s_\pi}({\frak n}_+)$ module by 
$u\otimes v y={\chi_h^{s_\pi}}(v)uy$, where $u\in Z_h^{s_\pi}({\frak g}),~~
v\in U_h^{s_\pi}({\frak n}_+),~~y\in Z_h^{s_\pi}({\frak g})/Z_{h,*}$. We denote 
this module by $(Z_h^{s_\pi}({\frak g})/Z_{h,*})_{\chi_h^{s_\pi}}$

The following result is another way of expressing Theorem ${\rm G_h}$. 
\vskip 0.3cm
\noindent
{\bf Theorem $\bf H_h$}
{\em Let $V_h$ be any $U_h^{s_\pi}({\frak g})$ module. Then $V_h$ is a Whittaker 
module if and only if one has an isomorphism
\begin{equation}\label{whittindh}
V_h\cong U_h^{s_\pi}({\frak g})\otimes_{Z_h^{s_\pi}({\frak g})\otimes 
U_h^{s_\pi}({\frak n}_+)}(Z_h^{s_\pi}({\frak g})/Z_{h,*})_{\chi_h^{s_\pi}}
\end{equation}
of $U_h^{s_\pi}({\frak g})$ modules. Furthermore in such a case the ideal 
$Z_{h,*}$ is unique and is given by 
$Z_{h,*}=Z_{h,V}^{s_\pi}({\frak g})$, where $Z_{h,V}^{s_\pi}({\frak g})$ is 
defined by (\ref{zvh}).}
\vskip 0.3cm

Note that the 
question of reducibility for $U_h^{s_\pi}({\frak g})$ modules which are free as 
${\Bbb C}[[h]]$ modules does not
make any sense since if $V_h$ is such a module then $hV_h$ is a proper 
subrepresentation. However it is natural to
study indecomposable modules for $U_h^{s_\pi}({\frak g})$. Below we describe all 
indecomposable Whittaker modules for 
$U_h^{s_\pi}({\frak g})$. First we determine all the Whittaker vectors in a 
Whittaker module for $U_h^{s_\pi}({\frak g})$.
\vskip 0.3cm
\noindent
{\bf Theorem $\bf K_h$}
{\em Let $V_h$ be any $U_h^{s_\pi}({\frak g})$ module with a cyclic Whittaker 
vector $w_h\in V_h$. Then any $v\in V_h$ is
a Whittaker vector if and only if $v$ is of the form $v=uw_h$, where $u\in 
Z_h^{s_\pi}({\frak g})$. Thus the 
space of all Whittaker vectors in $V_h$ is a cyclic $Z_h^{s_\pi}({\frak g})$ 
module which is isomorphic to 
$Z_h^{s_\pi}({\frak g})/Z_{h,V}^{s_\pi}({\frak g})$.}
\vskip 0.3cm

If $V_h$ is any $U_h^{s_\pi}({\frak g})$ module then ${\rm End}_{U_h}~V_h$ 
denotes the algebra of operators on $V_h$ which
commute with the action of $U_h^{s_\pi}({\frak g})$. If $\pi_V: 
U_h^{s_\pi}({\frak g})\rightarrow {\rm End}~V_h$ is the representation 
defining the $U_h^{s_\pi}({\frak g})$ module structure on $V_h$ then clearly 
$\pi_V(Z_h^{s_\pi}({\frak g}))\subseteq {\rm End}_{U_h}~V_h$.
Furthermore it is also clear that $\pi_V(Z_h^{s_\pi}({\frak g}))\cong 
Z_h^{s_\pi}({\frak g})/Z_{h,V}^{s_\pi}({\frak g})$.
\vskip 0.3cm
\noindent
{\bf Theorem $\bf L_h$}
{\em Assume that $V_h$ is a Whittaker module. Then ${\rm 
End}_{U_h}~V_h=\pi_V(Z_h^{s_\pi}({\frak g}))$. In particular one has an
isomorphism
$$
{\rm End}_{U_h}~V_h\cong Z_h^{s_\pi}({\frak g})/Z_{h,V}^{s_\pi}({\frak g}).
$$
Note that ${\rm End}_{U_h}~V_h$ is commutative.}
\vskip 0.3cm

Now one can describe all indecomposable Whittaker modules for 
$U_h^{s_\pi}({\frak g})$. A homomorphism
$$
{\xi_h}: Z_h^{s_\pi}({\frak g})\rightarrow {\Bbb C}[[h]]
$$
is called a central character. Given a central character ${\xi_h}$ let 
$Z_{h,{\xi_h}}^{s_\pi}({\frak g})={\rm Ker}~{\xi_h}$ so that
$Z_{h,{\xi_h}}^{s_\pi}({\frak g})$ is a typical central ideal in 
$Z_h^{s_\pi}({\frak g})$.

If $V_h$ is any $U_h^{s_\pi}({\frak g})$ module one says that $V_h$ admits an 
infinitesimal character, and ${\xi_h}$ is its
infinitesimal character, if ${\xi_h}$ is a central character such that 
$uv={\xi_h}(u)v$ for all $u\in Z_h^{s_\pi}({\frak g}),~
v\in V_h$. 

Given a central character ${\xi_h}$ let ${\Bbb 
C}[[h]]_{{\xi_h},{\chi_h^{s_\pi}}}$ be the 1-dimensional 
$Z_h^{s_\pi}({\frak g})\otimes U_h^{s_\pi}({\frak n}_+)$ module defined so that 
if $u\in Z_h^{s_\pi}({\frak g}),~v\in U_h^{s_\pi}({\frak n}_+),~
y\in {\Bbb C}[[h]]_{{\xi_h},{\chi_h^{s_\pi}}}$ then $u\otimes v 
y={\xi_h}(u){\chi_h^{s_\pi}}(v)y$. Also let
$$
Y_{{\xi_h},{\chi_h^{s_\pi}}}=U_h^{s_\pi}({\frak g})\otimes_{Z_h^{s_\pi}({\frak 
g})\otimes U_h^{s_\pi}({\frak n}_+)}{\Bbb C}[[h]]_{{\xi_h},{\chi_h^{s_\pi}}} .
$$
It is clear that $Y_{{\xi_h},{\chi_h^{s_\pi}}}$ admits an infinitesimal 
character and ${\xi_h}$ is that character.
\vskip 0.3cm
\noindent
{\bf Theorem $\bf M_h$}
{\em Let $V_h$ be any Whittaker module for $U_h^{s_\pi}({\frak g})$. Then the 
following conditions are equivalent:

(1) $V_h$ is an indecomposable $U_h^{s_\pi}({\frak g})$ module.

(2) $V_h$ admits an infinitesimal character.

(3) The corresponding ideal given by Theorem $\rm G_h$ is a maximal ideal.

(4) The space of Whittaker vectors in $V_h$ is 1--dimensional (over ${\Bbb 
C}[[h]]$).

(5) The centralizer ${\rm End}_{U_h}~V_h$ reduces to ${\Bbb C}[[h]]$.

(6) $V_h$ is isomorphic to $Y_{{\xi_h},{\chi_h^{s_\pi}}}$ for some central 
character ${\xi_h}$.}
\vskip 0.3cm
\noindent
{\em Proof.} The equivalence of (2), (3), (4) and (5) follows from Theorems $\rm 
K_h$ and $\rm L_h$. One gets the equivalence with 
(7) by Theorem $\rm H_h$. It remains to relate (2)--(6) with (1).
 
First we shall prove that (1) implies (4). Suppose that $V_h$ is indecomposable. 
Then the corresponding Whittaker module $V=V_h/hV_h$ for $U({\frak g})$ is
irreducible, and hence by Theorem M the space of Whittaker vectors of $V$ is 
one--dimensional. We denote this space by $Wh(V)$. Let $w_h$ be a 
cyclic Whittaker vector for $V_h$. Denote $Wh'(V_h)={\Bbb C}[[h]]w_h$. Then 
clearly $Wh'(V_h)\mbox{ mod h }=Wh(V)$. On the other hand
if $Wh(V_h)$ is the space of all Whittaker vectors in $V_h$ then $Wh(V_h)\mbox{ 
mod h }=Wh(V)$. Now from Lemma \ref{mainl}
applied to $A=Wh(V_h),~B=Wh'(V_h)$ it follows that $Wh(V_h)=Wh'(V_h)={\Bbb 
C}[[h]]w_h$.

Now we shall prove that (6) implies (1). Assume that (6) is satisfied. Then by 
construction the 
corresponding Whittaker module $V=V_h/hV_h$ for $U({\frak g})$ is isomorphic to 
$Y_{\xi,\chi}$ for some 
central character $\xi$. Now suppose that $V_h$ is decomposable. Then $V$ must 
be reducible which is impossible by
Theorem M. This completes the proof of Theorem $\rm M_h$.


\subsection{Coxeter realizations of quantum groups and Drinfeld twist}

\setcounter{equation}{0}
\setcounter{theorem}{0}

In this section we show that the Coxeter realizations $U_h^{s_\pi}({\frak g})$ 
of the quantum group $U_h({\frak g})$
are connected with quantizations of some nonstandard bialgebra structures on 
$\frak g$. At the quantum level 
changing bialgebra structure corresponds to the so--called Drinfeld twist. We 
shall consider a particular class
of such twists described in the following proposition.
\begin{proposition}{\bf (\cite{ChP}, Proposition 4.2.13)}\label{twdef}
Let $(A,\mu , \imath , \Delta , \varepsilon , S)$ be a Hopf algebra over a 
commutative ring. Let ${\cal F}$ be an invertible element of $A\otimes A$
such that 
\begin{equation}\label{twist}
\begin{array}{l}
{{\cal F}}_{12}(\Delta \otimes id)({{\cal F}})={{\cal F}}_{23}(id \otimes 
\Delta)({{\cal F}}),\\
\\
(\varepsilon \otimes id)({{\cal F}})=(id \otimes \varepsilon )({{\cal F}})=1.
\end{array}
\end{equation}
Then, $v=\mu (id\otimes S)({{\cal F}})$ is an invertible element of $A$ with
$$
v^{-1}=\mu (S\otimes id)({{\cal F}}^{-1}).
$$

Moreover , if we define $\Delta^{{\cal F}}:A\rightarrow A\otimes A$ and 
$S^{{\cal F}}:A\rightarrow A$ by
$$
\Delta^{{\cal F}}(a)={{\cal F}}\Delta(a){{\cal F}}^{-1},~~S^{{\cal 
F}}(a)=vS(a)v^{-1},
$$
then $(A,\mu , \imath , \Delta^{{\cal F}} , \varepsilon , S^{{\cal F}})$ is a 
Hopf algebra denoted by $A^{{\cal F}}$
and called the twist of $A$ by ${{\cal F}}$.
\end{proposition}

\begin{corollary}{\bf (\cite{ChP}, Corollary 4.2.15)}
Suppose that $A$ and ${{\cal F}}$ as in Proposition \ref{twdef}, but assume in 
addition that $A$ is quasitriangular
with universal R--matrix ${\cal R}$. Then $A^{{\cal F}}$ is quasitriangular with 
universal R--matrix
\begin{equation}\label{rf}
{\cal R}^{{\cal F}}={{\cal F}}_{21}{\cal R}{{\cal F}}^{-1},
\end{equation}
where ${{\cal F}}_{21}=\sigma {{\cal F}}$.
\end{corollary}

Fix a Coxeter element $s_\pi\in W$, $s_\pi=s_{\pi (1)}\ldots s_{\pi (l)}$.
Consider the twist of the Hopf algebra $U_h({\frak g})$ by the element
\begin{equation}\label{Ftw}
{{\cal F}}=exp(-h\sum_{i,j=1}^l \frac{n_{ji}}{d_j}Y_i\otimes Y_j) \in U_h({\frak 
h})\otimes U_h({\frak h}),
\end{equation}
where $n_{ij}$ is a solution of the corresponding equation (\ref{eqpi}).

This element satisfies conditions (\ref{twist}), and so $U_h({\frak g})^{{\cal 
F}}$ is a quasitriangular 
Hopf algebra with the universal R--matrix ${\cal R}^{{\cal F}}={{\cal 
F}}_{21}{\cal R}{{\cal F}}^{-1}$, 
where ${\cal R}$ is given by (\ref{univr}). We shall explicitly calculate the 
element ${\cal R}^{{\cal F}}$.
Substituting (\ref{univr}) and (\ref{Ftw}) into (\ref{rf}) and using 
(\ref{roots-cart}) we obtain
$$
\begin{array}{l}
{\cal R}^{{\cal F}}=exp\left[ h(\sum_{i=1}^l(Y_i\otimes H_i)+
\sum_{i,j=1}^l (-\frac{n_{ij}}{d_i}+\frac{n_{ji}}{d_j})Y_i\otimes Y_j) 
\right]\times \\
\prod_{\beta}
exp_{q_{\beta}^{-1}}[(q-q^{-1})a(\beta)^{-1}X_{\beta}^+e^{hK\beta^\vee} \otimes 
e^{-hK^*\beta^\vee}X_{\beta}^-],
\end{array}
$$
where $K$ is defined by (\ref{Kdef}).

Equip $U_h^{s_\pi}({\frak g})$ with the comultiplication given by :
$\Delta_{s_\pi}(x)=(\psi_{\{ n\}}^{-1}\otimes \psi_{\{ n\}}^{-1})\Delta_h^{{\cal 
F}}(\psi_{\{ n\}}(x))$.
Then $U_h^{s_\pi}({\frak g})$ becomes a quasitriangular Hopf algebra with the 
universal R--matrix
${\cal R}^{s_\pi}=\psi_{\{ n\}}^{-1}\otimes \psi_{\{ n\}}^{-1}{\cal R}^{{\cal 
F}}$. Using equation
(\ref{eqpi}) and Lemma \ref{tmatrel} this R--matrix may be written as follows
\begin{equation}\label{rmatrspi}
\begin{array}{l}
{\cal R}^{s_\pi}=exp\left[ h(\sum_{i=1}^l(Y_i\otimes H_i)+
\sum_{i=1}^l \frac{1+s_\pi}{1-s_\pi }H_i\otimes Y_i) \right]\times \\
\prod_{\beta}
exp_{q_{\beta}^{-1}}[(q-q^{-1})a(\beta)^{-1}e_{\beta} \otimes 
e^{h\frac{1+s_\pi}{1-s_\pi} \beta^\vee}f_{\beta}].
\end{array}
\end{equation}

The element ${\cal R}^{s_\pi}$ may be also represented in the form
\begin{equation}\label{rmatrspi'}
\begin{array}{l}
{\cal R}^{s_\pi}=exp\left[ h(\sum_{i=1}^l(Y_i\otimes H_i)\right]\times \\
\prod_{\beta}
exp_{q_{\beta}^{-1}}[(q-q^{-1})a(\beta)^{-1}e_{\beta}e^{-h\frac{1+s_\pi}{1-s_\pi
}\beta^\vee}\otimes f_\beta]
exp\left[ h(\sum_{i=1}^l \frac{1+s_\pi}{1-s_\pi }H_i\otimes Y_i)\right] .
\end{array}
\end{equation}

The comultiplication $\Delta_{s_\pi}$ is given on generators by
$$
\begin{array}{l}
\Delta_{s_\pi}(H_i)=H_i\otimes 1+1\otimes H_i,\\
\\
\Delta_{s_\pi}(e_i)=e_i\otimes e^{hd_i\frac{2}{1-s_\pi}H_i}+1\otimes e_i,\\
\\
\Delta_{s_\pi}(f_i)=f_i\otimes 
e^{-hd_i\frac{1+s_\pi}{1-s_\pi}H_i}+e^{-hd_iH_i}\otimes f_i.
\end{array}
$$

Note that the Hopf algebra $U_h^{s_\pi}({\frak g})$ is a quantization of the 
bialgebra structure on $\frak g$
defined by the cocycle
\begin{equation}\label{cocycles}
\delta (x)=({\rm ad}_x\otimes 1+1\otimes {\rm ad}_x)2r^{s_\pi}_+,~~ 
r^{s_\pi}_+\in {\frak g}\otimes {\frak g},
\end{equation}
where $r^{s_\pi}_+=r_+ + \frac 12 \sum_{i=1}^l \frac{1+s_\pi}{1-s_\pi 
}H_i\otimes Y_i$, and $r_+$ is given by (\ref{rcl}).
 
We shall also need the following property of the antipode $S^{s_\pi}$ of 
$U_h^{s_\pi}({\frak g})$.
\begin{proposition}\label{sqant}
The square of the antipode $S^{s_\pi}$ is an inner automorphism of 
$U_h^{s_\pi}({\frak g})$ given by
$$
(S^{s_\pi})^2(x)=e^{2h\rho^\vee}xe^{-2h\rho^\vee},
$$
where $\rho^\vee=\sum_{i=1}^lY_i$.
\end{proposition}
{\em Proof.} 
First observe that by Proposition \ref{twdef} the antipode of 
$U_h^{s_\pi}({\frak g})$ has the form:
$S^{s_\pi}(x)=\psi_{\{ n\}}^{-1}(vS_h(\psi_{\{ n\}}(x))v^{-1})$, where 
$$
v=exp(h\sum_{i,j=1}^l \frac{n_{ji}}{d_j}Y_iY_j).
$$
Therefore $(S^{s_\pi})^2(x)=\psi_{\{ n\}}^{-1}(vS_h(v^{-1})S_h^2(\psi_{\{ 
n\}}(x))S_h(v)v^{-1})$.
Note that $S_h(v)=v$, and hence $(S^{s_\pi})^2(x)=\psi_{\{ 
n\}}^{-1}(S_h^2(\psi_{\{ n\}}(x)))$.

Finally observe that from explicit formulas for the antipode of $U_h({\frak g})$ 
it follows that 
$S_h^2(x)=e^{2h\rho^\vee}xe^{-2h\rho^\vee}$. This completes the proof.

In conclusion we note that using Corollary \ref{rep} and the isomorphism 
$\psi_{\{ n\}}$ one can define finite--dimensional representations of 
$U_h^{s_\pi}({\frak g})$.


\subsection{Quantum deformation of the Toda lattice}\label{toda}

\setcounter{equation}{0}
\setcounter{theorem}{0}

Recall that one of
the main applications of the algebra $W({\frak b}_-)$ is the quantum Toda 
lattice \cite{K'}.
Let $\overline \chi : {\frak n}_- \rightarrow {\Bbb C}$ be a non--singular 
character of
the opposite nilpotent subalgebra ${\frak n}_-$. We denote the character of 
$N_-$ corresponding to 
$\overline \chi$ by the same letter. The algebra $U({\frak b}_-)$ naturally acts 
by differential
operators in the space $C^\infty ({\Bbb C}_{\overline \chi}\otimes_{N_-}{B_-})$. 
This space may be
identified with $C^\infty (H)$.
Let $D_1,\ldots ,D_l$ be the differential operators on $C^\infty (H)$ which 
correspond to the elements 
$I_1^\chi ,\ldots , I_l^\chi\in W({\frak b}_-)$. Denote by $\varphi$ the 
operator of multiplication in 
$C^\infty (H)$ by the function $\varphi (e^h)=e^{\rho(h)}$, where $h\in {\frak 
h}$. The operators $M_i=\varphi D_i\varphi^{-1}, i=1,\ldots l$
are called the quantum Toda Hamiltonians. Clearly, they commute with each other.

In particular if $I$ is the quadratic Casimir element then the corresponding 
operator $M$ is 
the well--known second--order differential operator:
$$
M=\sum_{i=1}^l \partial_i^2
+\sum_{i=1}^l \chi(X_{\alpha_i})\overline 
\chi(X_{-\alpha_i})e^{-\alpha_i(h)}+(\rho,\rho),
$$
where $\partial_i=\frac{\partial}{\partial y_i}$, and $y_i,~i=1,\ldots l$ is an 
ortonormal basis of $\frak h$.

Using the algebra $W_h({\frak b}_-)$ we shall construct quantum group analogues 
of the Toda Hamiltonians.
A slightly different approach has been recently proposed in \cite{Et}.

Denote by $A$ the space of linear functions on 
${\Bbb C}[[h]]_{\overline \chi_h^{s_\pi}}\otimes_{U_h^{s_\pi}({\frak 
n}_-)}U_h^{s_\pi}({\frak b}_-)$, where
${\Bbb C}[[h]]_{\overline \chi_h^{s_\pi}}$ is the one--dimensional 
$U_h^{s_\pi}({\frak n}_-)$ module
defined by ${\overline \chi_h^{s_\pi}}$.
Note that 
${\Bbb C}[[h]]_{\overline \chi_h^{s_\pi}}\otimes_{U_h^{s_\pi}({\frak 
n}_-)}U_h^{s_\pi}({\frak b}_-)\cong U_h^{s_\pi}({\frak h})$
as a linear space. Therefore $A=U_h^{s_\pi}({\frak h})^*$.
The algebra $U_h^{s_\pi}({\frak b}_-)$ naturally acts on 
${\Bbb C}[[h]]_{\overline \chi_h^{s_\pi}}\otimes_{U_h^{s_\pi}({\frak 
n}_-)}U_h^{s_\pi}({\frak b}_-)$ 
by multiplications from the right. This action induces an $U_h^{s_\pi}({\frak 
b}_-)$--action in the space $A$.
We denote this action by $L$, $L:U_h^{s_\pi}({\frak b}_-)\rightarrow {\rm 
End}A$. Clearly, this action generates 
an action of the algebra $W_h({\frak b}_-)$ on $A$.

To construct deformed Toda Hamiltonians we shall use certain elements in  
$W_h({\frak b}_-)$.
These elements may be described as follows.
Let $\mu : U_h^{s_\pi}({\frak g}) \rightarrow {\Bbb C}[[h]]$ be a map such that 
$\mu(uv)=\mu(vu)$. By Proposition 
\ref{sqant} $(S^{s_\pi})^2(x)=e^{2h\rho^\vee}xe^{-2h\rho^\vee}$. Hence from 
Remark 1 in \cite{D} (see also \cite{FRT}) it follows that  
$(id\otimes \mu)({\cal R}_{21}^{s_\pi}{\cal R}^{s_\pi}(1\otimes 
e^{2h\rho^\vee}))$, where
${\cal R}_{21}^{s_\pi}=\sigma {\cal R}^{s_\pi}$, is a central element.
In particular, for any finite--dimensional $\frak g$--module $V$ the element
\begin{equation}\label{centrelv}
C_V=(id\otimes tr_V)({\cal R}_{21}^{s_\pi}{\cal R}^{s_\pi}(1\otimes 
e^{2h\rho^\vee})),
\end{equation}
where $tr_V$ is the 
trace in $V[[h]]$, is central in $U_h^{s_\pi}({\frak g})$.

Using formulas (\ref{rmatrspi}) and (\ref{rmatrspi'}) we can easily compute 
elements 
$\rho_{\chi_h^{s_\pi}}(C_V)\in W_h({\frak b}_-)$.
For every finite--dimensional $\frak g$--module $V$ we have
\begin{equation}\label{todah}
\begin{array}{l}
\rho_{\chi_h^{s_\pi}}(C_V)=(id\otimes tr_V)( e^{t_0}\prod_{\beta}
exp_{q_{\beta}^{-1}}[(q-q^{-1})a(\beta)^{-1}f_\beta \otimes 
e_{\beta}e^{-h\frac{1+s_\pi}{1-s_\pi}\beta^\vee}]\times \\
\\
e^{t_0}\prod_{\beta}
exp_{q_{\beta}^{-1}}[(q-q^{-1})a(\beta)^{-1}{\chi_h^{s_\pi}}(e_{\beta}) \otimes 
e^{h\frac{1+s_\pi}{1-s_\pi} \beta^\vee}f_{\beta}](1\otimes e^{2h\rho^\vee})),
\end{array}
\end{equation}
where $t_0=h\sum_{i=1}^l(Y_i\otimes H_i)$.

We denote by $W_h^{Rep}({\frak b}_-)$ the subalgebra in $W_h({\frak b}_-)$ 
generated by the elements
$\rho_{\chi_h^{s_\pi}}(C_V)$, where $V$ runs through all finite--dimensional 
representations of $\frak g$.
Note that for every finite--dimensional $\frak g$--module $V$ 
$\rho_{\chi_h^{s_\pi}}(C_V)$ is a 
polynomial in noncommutative elements $f_i,~e^{hx},~x\in {\frak h}$.

Now we shall realize elements of $W_h^{Rep}({\frak b}_-)$ as difference 
operators.
Let $H_h\in U_h^{s_\pi}({\frak h})$ be the subgroup generated by elements 
$e^{hx},~x\in {\frak h}$.
A difference operator on $A$ is an operator $T$ of the form
$T=\sum f_iT_{x_i}$ (a finite sum), where $f_i \in A$, and for every $y\in H_h~$ 
$T_{x}f(y)=(ye^{hx}),~x\in {\frak h}$.
\begin{proposition}{\bf (\cite{Et}, Proposition 3.2)}\label{diffh}
For any $Y \in U_h^{s_\pi}({\frak b}_-)$, 
which is a polynomial in noncommutative elements $f_i,~e^{hx},~x\in {\frak h}$,
the operator $L(Y)$ is a difference operator on $A$.
In particular, the operators $L(I),~I\in W_h^{Rep}({\frak b}_-)$ are mutually 
commuting 
difference operators on $A$.
\end{proposition}
{\em Proof.} It suffices to verify that $L(f_i)$ are difference operators on 
$H_h$.
Indeed,
$$
L(f_i)f(e^{hx})=f(e^{hx}f_i)=e^{-h\alpha_i(x)}f(f_ie^{hx})=\overline 
\chi_h^{s_\pi}(f_i)e^{-h\alpha_i(x)}f(e^{hx}).
$$
This completes the proof.

Let $\jmath : H_h\rightarrow U_h^{s_\pi}({\frak h})$ be the canonical embedding. 
Denote $A_h=\jmath^*(A)$.
Let $T$ be a difference operator on $A$. Then one can define a difference 
operator $\jmath^*(T)$ on the space
$A_h$ by $\jmath^*(T)f(y)=T(\jmath(y))$.  

Let $D_i^h=\jmath^*(L(\rho_{\chi_h^{s_\pi}}(C_{V_i})))$, where $V_i,~i=1,\ldots 
l$ are the fundamental representations of
$\frak g$. 
Denote by $\varphi_h$ the operator of multiplication in 
$A_h$ by the function $\varphi_h (e^{hx})=e^{h\rho(x)}$, where $x\in {\frak h}$. 
The operators $M_i^h=\varphi_h D_i^h\varphi^{-1}_h, i=1,\ldots l$
are called the quantum deformed Toda Hamiltonians.

From now on we suppose that $\pi=id$ and that the ordering of positive roots 
$\Delta_+$ is fixed as in 
Proposition \ref{rootsh}. We denote $s_{id}=s$.
Now using formula (\ref{todah}) we outline computation of the operators $M_i^h$. 
This computation is simplified by
the following lemma.
\begin{lemma}{ \bf (\cite{Et}, Lemma 5.2)} 
Let $X=f_{\gamma_1}...f_{\gamma_n}$. If the roots
$\gamma_1,...,\gamma_n$ are not all simple 
then $L(X)=0$. 
Otherwise, if $\gamma_i=\alpha_{k_i}$, then 
$$
\jmath^*(L(X))f(e^{hy})=e^{-h(\sum\alpha_{k_i},y)}f(e^{hy})\prod_i\overline 
\chi_h^{s}(f_{k_i})
$$ 
\end{lemma}
{\em Proof }  follows immediately from Proposition \ref{rootsh} and the 
arguments used in the proof of 
Proposition \ref{diffh}. 

Using this lemma we obtain that if $\beta$ is not a simple root then the term in 
(\ref{todah}) containing 
root vector $f_\beta$ gives a trivial contribution to the operators
$L(\rho_{\chi_h^{s}}(C_{V_i}))$. Note also that by Proposition \ref{rootsh} 
${\chi_h^{s}}(e_\beta)=0$ if 
$\beta$ is not a simple root. Therefore from formula (\ref{todah}) we have
\begin{equation}
\begin{array}{l}
L(\rho_{\chi_h^{s}}(C_{V_i}))=\\
\\
L(id\otimes tr_V)( e^{t_0}\prod_{i}
exp_{q^{-2d_i}}[(q_i-q_i^{-1})f_i \otimes e_ie^{-hd_i\frac{1+s}{1-s}H_i}]\times 
\\
\\
e^{t_0}\prod_{i}
exp_{q^{-2d_i}}[(q_i-q_i^{-1}){\chi_h^{s}}(e_i) \otimes 
e^{hd_i\frac{1+s}{1-s}H_i}f_i](1\otimes e^{2h\rho^\vee})).
\end{array}
\end{equation}

In particular, let ${\frak g}=sl(n)$, $V_1=V$ the fundamental representation of 
$sl(n)$. 
Then direct calculation gives
$$
M_1f(e^{hy})=\left( \sum_{j=1}^n T_j^2-
(q-q^{-1})^2\sum_{i=1}^{n-1}{\chi_h^{s}}(e_i){\overline \chi_h^{s}}(f_i)
e^{-h(y,\alpha_i)}T_{i+1}T_i\right) f(e^{hy}),
$$
where $T_i=T_{\omega_i}$, ${\omega_i}$ are the weights of $V$. 
The last expression coincides with formula (5.7) obtained in \cite{Et}.


\end{document}